\documentclass[preprint,3p,12pt]{elsarticle}
\journal{Information Sciences}

\usepackage{amsmath,comment,amsfonts}
\usepackage{bm}
\usepackage{array}
\usepackage{amssymb}
\usepackage{amscd}
\usepackage{multirow}
\usepackage{calrsfs}
\usepackage{url}
\usepackage{color}
\usepackage{soul}
\usepackage{xcolor,yfonts, graphicx,csquotes,amsfonts,amsmath,amssymb,bbm,mathrsfs,bbm}
\usepackage[]{algorithm2e}
\usepackage{csquotes}
\usepackage[english]{babel}
\everymath{\displaystyle}
\usepackage[normalem]{ulem}
\usepackage{float,dirtytalk}
\usepackage{tasks}
\def\ml{l\kern-0.035cm\char39\kern-0.03cm}

\newtheorem {definition}{Definition}[section]

\newtheorem {theorem}[definition]{Theorem}
\newtheorem {lema}[definition]{Lemma}

\newtheorem {proposition}[definition]{Proposition}

\newtheorem {example}[definition]{Example}
\newtheorem {remark}[definition]{Remark}
\newtheorem {corollary}[definition]{Corollary}

\newcommand{\proof}{\medskip \noindent {\bf Proof. \ \ }}
\newcommand{\collection}{\mathcal{D}}
\newcommand{\collectionGSF}{\mathcal{E}}

\newcommand{\size}{\mathsf{s}}

\newcommand{\counting}{{\#}}

\newcommand{\x}{\mathbf{x}}
\newcommand{\y}{\mathbf{y}}
\newcommand{\z}{\mathbf{z}}
\newcommand{\w}{\mathbf{w}}

\newcommand{\bin}{\mathbf{1}}

\newcommand\sA{\mathsf{A}}
\newcommand\aA[2][E]{\sA(#2|#1)}%conditional aggregation
\newcommand\aAi[3][E]{\sA^{#3}(#2|#1)}

\newcommand\cA{{\mathcal A}}

\newcommand\bA{{\mathscr {A}}}

\newcommand\cN{{{\mathcal N}_m}}
\newcommand\cE{{\mathcal{E}}}

\newcommand\bM{{\mathbf M}}
\newcommand\Ai[1]{E_{(#1)}}

\def\lm{\mu(\{\x>\alpha\})}

\def\gsf{\mu_{\bA}(\x,\alpha)}

\def\gsf#1{\mu_{\mathcal{A}^{#1}}(\x,\alpha)}

\def\gsf#1#2#3{\mu_{\mathcal{A}^{#1}}(#2,#3)}

\usepackage{xcolor}
\definecolor{darkgreen}{rgb}{0,1,0.1}

\definecolor{brown}{rgb}{0.82,0.41,0.15}%chocolate

\definecolor{brown}{rgb}{0.43, 0.21, 0.1}

\usepackage{enumerate}
\newcommand\gout{\bgroup\markoverwith{\textcolor{green}{\rule[0.5ex]{5pt}{1.5pt}}}\ULon}

\title{}

\date{ }
\usepackage{ulem}%skrtnutie textu

\begin{document}
%\maketitle

\begin{frontmatter}
\title{Survival functions versus conditional aggregation-based survival functions on~discrete~space}
\author{{Basarik Stanislav}\fnref{fn1}}\ead{stanislav.basarik@student.upjs.sk} 

\author{{Borzová Jana}\fnref{fn1}}\ead{jana.borzova@upjs.sk}

\author{{Halčinová Lenka}\corref{cor1}\fnref{fn1}}\ead{lenka.halcinova@upjs.sk}\cortext[cor1]{lenka.halcinova@upjs.sk}

\address[fn1]{Institute of Mathematics,
P.~J. \v{S}af\'arik University in Ko\v sice,
Jesenn\'a 5, 040 01 Ko\v{s}ice, Slovakia}
\fntext[fn2]{Supported by the grants APVV-16-0337, VEGA 1/0657/22, bilateral call Slovak-Poland
grant scheme No. SK-PL-18-0032 and grant scheme
VVGS-PF-2021-1782.
}

\begin{abstract}
In this paper we deal with conditional aggregation-based survival functions recently introduced by Boczek et al. (2020). The concept is worth to study because of its 
possible implementation in real-life situations and mathematical theory as well. 
The aim of this paper is the comparison of this new notion with the standard survival function.
We state sufficient and necessary conditions under which the generalized and the standard survival function equal.
The main result is the characterization of the family of conditional aggregation operators (on discrete space) for which these functions coincide. 

\end{abstract}

\end{frontmatter}

\section{Introduction}

We continue to study the~novel survival functions introduced in~\cite{BoczekHalcinovaHutnikKaluszka2020} as a~generalization of size-based level measure developed for the~use in nonadditive analysis in~\cite{BorzovaHalcinovaSupina2018,Halcinova2017,HalcinovaHutnikKiselakSupina2019}.  The~concept appeared initially in time-frequency analysis~\cite{DoThiele2015}. As the main result, in Theorem~\ref{thm:characterization2} we show that the~generalized survival function is equal to the original notion (for any monotone measure and any input vector) just in very particular case. The concept of the novel survival function is useful in many real-life situations and pure theory as well. 
In fact, the~standard survival function (also known in the literature as the standard level measure~\cite{HalcinovaHutnikKiselakSupina2019}, strict level measure~\cite{BorzovaHalcinovaHutnik2020} or decumulative distribution
function~\cite{grabisch2016set})
is the~crucial ingredient of many definitions in mathematical analysis. Many well-known integrals are based on the~survival function, e.g. the~Choquet integral, the~Sugeno integral, the~Shilkret integral, the~seminormed integral~\cite{BorzovaHalcinovaHutnik2020}, universal integrals~\cite{KlementMesiarPap2010}, etc. Also, the convergence of a~sequence of functions in measure is based on the~same concept. Hence a~reasonable ge\-ne\-ra\-li\-za\-tion of~the~survival function leads to the generalizations of all mentioned concepts. For more on applications of the~generalized survival function, see~\cite{BoczekHalcinovaHutnikKaluszka2020,DoThiele2015}. 
\bigskip

Due to the~number of factors needed in the~definition of the~generalized survival function, it is quite difficult to understand this concept. In order to understand it more deeply, in~the~following we shall focus on the graphical visualization of inputs, see~\cite{BorzovaHalcinovaSupina2021}. Moreover, the graphical representation will help us to formulate basic results of this paper.
In~the~whole paper, we restrict ourselves to discrete settings. We consider finite basic set $$[n]:=\{1,2,\dots, n\}, \,\, n\geq1$$ and a~monotone measure~$\mu$ on $2^{[n]}$.  If $\x=(x_1,\dots, x_n)$ is a~nonnegative real-valued function on  $[n]$, i.e., a~vector, then the \textit{survival function} (or standard survival function) of~the~vector $\x$ with respect to $\mu$, see~\cite{BoczekHalcinovaHutnikKaluszka2020,DuranteSempi2015}, is 
defined 
by 
$$\lm:=\mu\left(\{i\in [n]: x_i>\alpha\}\right), \quad \alpha\in[0,\infty).$$
For the~thorough exposition see Preliminaries. 
To avoid too abstract setting in the following visual representations, let us consider the~input vector $\x=(2,3,4)$ and the~monotone measure $\mu$ on~$ 2^{[3]}$ defined in~Table~\ref{tabulka_0}.
\par

\begin{table}
\renewcommand*{\arraystretch}{1.5}
\begin{center}
\begin{tabular}{|>{\centering\arraybackslash}m{1.25cm}|>{\centering\arraybackslash}m{1.5cm}|>{\centering\arraybackslash}m{1cm}|>{\centering\arraybackslash}m{1cm}|>{\centering\arraybackslash}m{1cm}|>{\centering\arraybackslash}m{1cm}|>{\centering\arraybackslash}m{1cm}|>{\centering\arraybackslash}m{1cm}|>{\centering\arraybackslash}m{1.5cm}|}
\hline
$E$ & $\{1,2,3\}$ & $\{2,3\}$ & $\{1,3\}$ & $\{1,2\}$ & $\{3\}$ & $\{2\}$ & $\{1\}$ & $\emptyset$ \\ \hline
$E^c$ & $\emptyset$ & $\{1\}$ & $\{2\}$ & $\{3\}$ & $\{1,2\}$ & $\{1,3\}$ & $\{2,3\}$ & $\{1,2,3\}$ \\ \hline
$\mu(E^c)$ & $0$ & $0.25$ & $0.25$ & $0.4$ & $0.75$ & $0.75$ & $0.75$ & $1$\\ \hline
$\max_{i\in E}x_i$ & $4$  & $4$ & $4$ & $3$ & $4$ & $3$ & $2$ & $0$\\\hline
$\sum_{i\in E}x_i$ & $9$  & $7$ & $6$ & $5$ & $4$ & $3$ & $2$ & $0$\\\hline
\end{tabular}\caption{Sample measure~$\mu$ and two conditional aggregation operators for vector $\x=(2,3,4)$}\label{tabulka_0}
\end{center}
\end{table}

\paragraph{\bf{The survival functions visual representation}}

We begin with a~nonstandard representation of standard survival function, as a~stepping stone to its~generalization. Before, let us introduce the following equivalent representation of survival function:
\begin{align}
\begin{split}\label{prepis2}
\lm=\mu([n]\setminus \{i\in [n]: x_i\leq \alpha\})&=\min\big\{\mu (E^c):(\forall i\in E)\,\, x_i\leq \alpha,\, E\in 2^{[n]}\big\}\\&=\min\big\{\mu (E^c): \max_{i\in E} x_i\leq \alpha,\, E\in 2^{[n]}\big\},
\end{split}
\end{align}
where $E^c=[n]\setminus E$, see motivation problem 1 in~\cite{BoczekHalcinovaHutnikKaluszka2020}. Let us start the visualization with inputs from Table~\ref{tabulka_0}.
\begin{figure}
\begin{center}
\begin{tabular}{m{7.8cm} m{6cm}}
{\vspace*{7pt}}\includegraphics[scale=1.1]{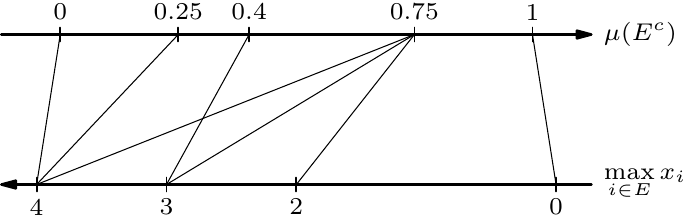} & \includegraphics[scale=1.1]{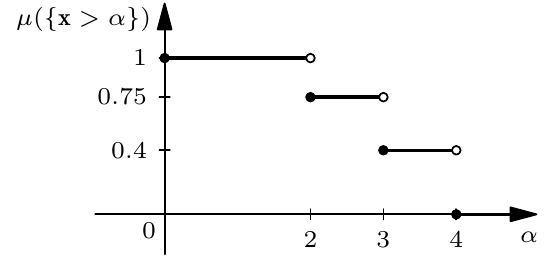}
\end{tabular}
\caption{The survival function visualization for $\x=(2,3,4)$ and $\mu$ given in~Table~\ref{tabulka_0}.}\label{MLM_vizualizacia}
\end{center}
\end{figure}
Let us depict all maximal values of~$\x$ on $E$, for each set $E\in2^{[3]}$ on the~lower axis, see left image of Figure~\ref{MLM_vizualizacia}, in decreasing order and the corresponding values of monotone measure of complement, i.e. $\mu(E^c)$,  on the~upper axis. 
In this picture of Figure~\ref{MLM_vizualizacia}, the~number on lower axis is linked with the~number on the~upper one via a~straight line once they correspond to the~same set, i.e., $a$ is linked with $b$ if there is $E\in2^{[3]}$ such that

$$a=\max\limits_{i\in E}x_i\hspace{0.5cm}\text{ and }\hspace{0.5cm} b=\mu(E^c).$$ 
Finally, the~value $\lm$ at some $\alpha\in[0,\infty)$ can be read from the left image of Figure~\ref{MLM_vizualizacia} considering the~minimal value on the~upper axis which is linked to a~value smaller than $\alpha$ (i.e., right-hand side value) on the~lower one. Thus considering an~illustrative example in the left image of Figure~\ref{MLM_vizualizacia}, the~value of survival function at $2{.}5$ is $0{.}75$. Indeed, there are just 2 values on the~right hand side of $2{.}5$, namely numbers 2 and 0. These are linked to $0{.}75$ and 1, respectively. Hence, $0{.}75$ is a~smaller one. The graph of survival function is in the right image of Figure~\ref{MLM_vizualizacia}.

\begin{figure}
\begin{center}
\begin{tabular}{m{7.2cm} m{7cm}}
\vspace*{7pt}
\includegraphics[scale=1.1]{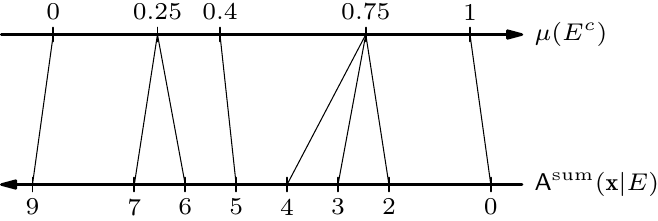} &
\includegraphics[scale=1.1]{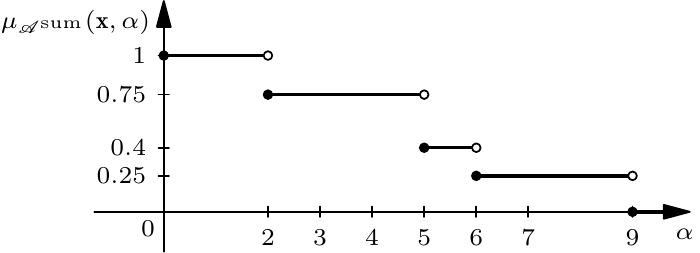}
\end{tabular}
\caption{Generalized survival function visualization for $\x=(2,3,4)$, $\mu$ given in~Table~\ref{tabulka_0} and $\sA=\sA^{\mathrm{sum}}$}\label{SLM_vizualizacia}
\end{center}
\end{figure}

\paragraph{\bf{The generalized survival functions visual representation}}
In the~modification of the survival function, the~previously described computational procedure stays. However, we allow to use any conditional aggregation operator, not just maximum operator. %on the~complement. 
The~standard example of conditional aggregation is the sum of components of~$\x$, 
see the last line in~Table~\ref{tabulka_0}
and the corresponding visualisation in Figure~\ref{SLM_vizualizacia}. 
Applying the~described computational procedure we obtain the~sum-based survival function of vector $\x$, i.e., the~generalized survival function  of vector $\x$ studied in~\cite{BoczekHalcinovaHutnikKaluszka2020,BorzovaHalcinovaSupina2018,Halcinova2017,HalcinovaHutnikKiselakSupina2019}. The formula linked to this procedure is the following:
\begin{eqnarray*}
\gsf{\mathrm{sum}}{\x}{\alpha}=\min\left\{\mu(E^c): \aAi[E]{\x}{\mathrm{sum}}\leq\alpha,\, E\in\collectionGSF\right\}\end{eqnarray*}
\noindent with $\aAi[E]{\x}{\mathrm{sum}}=\sum\limits_{i\in E}x_i$ and $\{\emptyset\}\subseteq\,\collectionGSF\subseteq 2^{[n]}$ (in the illustrative example $\collectionGSF=2^{[3]}$).

The~corresponding graph is: 
Considering discrete space, the~computation of the~generalized survival function studied in~\cite{BoczekHalcinovaHutnikKaluszka2020,BorzovaHalcinovaSupina2018,Halcinova2017,HalcinovaHutnikKiselakSupina2019} may be always represented via the~corresponding diagrams similar to those in~Figures~\ref{MLM_vizualizacia} and~\ref{SLM_vizualizacia}.  
\bigskip

Except for a better understanding of survival functions, the visual representation may help us to answer the problem of their indistinguishability.  With the introduction of novel survival function a natural question arises: When does the generalized survival function coincide with the survival function? 
The motivation for answering these questions is not only to know the relationship between mentioned concepts for given inputs, but it will help us to compare the corresponding integrals based on them, see~\cite[Definition 5.1, Definition 5.4]{BoczekHalcinovaHutnikKaluszka2020}. 
In the literature, there are known some families of conditional aggregation operators together with the collection $\collectionGSF$ when the generalized survival function equals to the survival function. In the following we list them:
\begin{itemize}
    \item  (cf.~\cite[Corollary 4.15]{HalcinovaHutnikKiselakSupina2019}) $\cA=\cA^{\rm{size}}$ with size~$\size$ being the weighted sum\footnote{ 
    ${\size}_{\counting,p}(\x)(E)= \left(\frac{1}{\counting(E)}\cdot\sum\limits_{x_i\in E}x_i^p\right)^{\frac{1}{p}} 
   $ \text{for}~$E\neq\emptyset$, ${\size}_{\counting,p}(\x)(\emptyset)=0$  and $p>0$.}, $\collection$ contains all singletons of~$[n]$ and $\collectionGSF=2^ {[n]}$;
        \item (cf.~~\cite[Example 4.2]{BoczekHalcinovaHutnikKaluszka2020} or~\cite[Section 5]{HalcinovaHutnikKiselakSupina2019}) $\cA=\cA^{\rm{max}}$ with $\collectionGSF=2^{[n]}$;
        \item (cf.~\cite[Proposition 4.6]{BoczekHalcinovaHutnikKaluszka2020}) $\cA=\cA^{\mu-\mathrm{ess}}$ with $\collectionGSF=2^{[n]}$.
\end{itemize}

Although the first two items appear to be different, in fact, under the above conditions, they are equal $\cA^{\rm{size}}=\cA^{\rm{max}}$. 
Settings of above mentioned examples lead to the survival function regardless of the choice of monotone measure $\mu$. However, the identity between generalized survival function and survival function may happen also for other 
families of conditional aggregation operators (a FCA for short), but with specific monotone measures, e.g. $\cA^{\mathrm{sum}}$ with the weakest monotone measure\footnote{$\mu_*\colon 2^ {[n]}\to[0,\infty)$ given by $$\mu_*(E)=\begin{cases}\mu([n]), & E=[n], \\ 0, & \textrm{otherwise}.\end{cases}$$} shrinks to survival function for any input vector $\x\in[0,\infty)^{[n]}$ and $\collectionGSF=2^{[n]}$. In this paper 
we shall treat the following problems:

\medskip

\noindent\textbf{Problem 1:} Let $\x\in[0,\infty)^{[n]}$, $\mu$ be a monotone measure on $2^{[n]}$, and $\cA$ be FCA. What are sufficient and necessary conditions on $\x$, $\mu$ and $\cA$ to hold $\gsf{}{\x}{\alpha}=\lm$?
\bigskip

\noindent\textbf{Problem 2:} Let $\x\in[0,\infty)^{[n]}$, and $\cA$ be FCA. What are sufficient and necessary conditions on $\x$ and $\cA$ to hold $\gsf{}{\x}{\alpha}=\lm$ for any monotone measure $\mu$?
\bigskip

\noindent\textbf{Problem 3:} Let $\cA$ be FCA. What are sufficient and necessary conditions on $\cA$ to hold $\gsf{}{\x}{\alpha}=\lm$ for any monotone measure $\mu$ and $\x\in[0,\infty)^{[n]}$?
\medskip

The~paper is organized as follows. We continue with preliminary section containing needed definitions and notations. In Section~3 we solve Problem 1, see e.g. Corollary~\ref{LM_SLM_coincide}, 
Corollary~\ref{corollary_1}, Remark~\ref{remark_i}, Proposition~\ref{summar} and Theorem~\ref{nutna_postacujuca_2}.  In Section~4 we provide quite~surprising result, see Theorem~\ref{thm:characterization2} that characterizes the family of conditional aggregation operators (in discrete setting) for which the~generalized survival function coincides with the~standard survival function. Thus we answer Problem 3. In Section~4 we also treat Problem~2, see Theorem~\ref{thm:characterization} and Theorem~\ref{thm:characterization_mon}.
Many our results are supported by appropriate examples.

\section{Background and preliminaries}

In order to be self-contained as far as possible, we recall in this section necessary definitions and all basic notations. In the whole paper, we restrict ourselves to discrete settings. As we have already mentioned, we shall consider a
finite set $$X=[n]:=\{1,2,\dots, n\},\,\, n\geq 1.$$
%\jb{ta 2 predtym bola kvoli tomu, ze v lit. sa to takto pouzivalo, zrejme preto, ze pre n=1 velmi nie je co riesit}
We shall denote by 
% $\Sigma$ the $\sigma$-algebra of $[n]$ and 
$2^{[n]}$ the power set of $[n]$.
A~{\it monotone} or \textit{nonadditive measure} on $2^{[n]}$  is a~nondecreasing set function $\mu\colon 2^{[n]}\to {{[0,\infty)}},$ i.e., $\mu(E)\le\mu(F)$ whenever $E\subseteq F$, with $\mu(\emptyset)=0.$ Moreover, we shall suppose $\mu([n])>0$. The set of monotone measures on $2^{[n]}$ we shall denote by $\mathbf{M}$.
The monotone measure satisfying the equality $\mu([n])=1$ will be called the \textit{normalized monotone measure} (also known as a capacity in~\cite{Weber1986}). In this paper we shall always work with monotone measures being defined on $2^{[n]}$, although, on several places the domain of $\mu$ can be smaller. Also, we shall need special properties of $\mu$ on a system $ \mathcal{S}\subseteq 2^{[n]}$. The monotone measure $\mu\in\mathbf{M}$ with the property $\mu(E)\neq\mu(F)$ for any $E,F\in \mathcal{S}\subseteq 2^{[n]}$, $E\neq F$  will be called \textit{strictly monotone measure} on $\mathcal{S}$.
The counting measure will be denoted by $\counting$. 
Further, we put $\max\emptyset=0$ and $\sum_{i\in\emptyset}x_i=0$.

We shall work with nonnegative real-valued vectors, we shall use the notation $\x=(x_1,\dots,x_n)$, $x_i\in [0,\infty)$, $i=1,2,\dots, n$.   The set $[0,\infty)^{[n]}$  is the family of all nonnegative real-valued functions on  $[n]$, i.e. vectors. For any  $\x=(x_1,\dots,x_n)\in[0,\infty)^{[n]}$ we denote by $(\cdot)$ a~permutation $(\cdot)\colon [n]\to[n]$ such that $x_{(1)}\leq x_{(2)}\leq\dots \leq x_{(n)}$ and $x_{(0)}=0$, $x_{(n+1)}=\infty$ by convention. Let us remark that the
permutation $(\cdot)$ need not be unique (this happens if there are some ties in the sample $(x_1,..., x_n)$, see~\cite{ChenMesiarLiStupnanova2017}). For a~fixed input vector $\x$ and a~fixed permutation $(\cdot)$  we shall denote by  $\Ai{i}$ the set of the form $\Ai{i}=\{(i),\dots,(n)\}$ for any $i\in [n]$ with the convention~$\Ai{n+1}=\emptyset$.
By $\bin_E$ we shall denote the indicator function of a set $E\subseteq Y$, $Y\subseteq [0,\infty)$,
i.e., $\bin_E(x)=1$ if $x\in E$ and $\bin_E(x)=0$ if $x\notin E$. Especially, $\bin_{\emptyset} (x)= 0$ for each $x\in Y$. 
We shall work with indicator function with respect to two different sets. We shall work with $Y=[n]$ when dealing with vectors (i.e. $\bin_E$ is a characteristic vector of $E\subseteq [n]$ in $\{0,1\}^{{[n]}}$) and $Y= [0,\infty)$ when dealing with survival functions.

In the following we list several definitions (adopted to discrete settings).
Firstly, the concept of the \textit{conditional aggregation operator} is presented. Its crucial feature is that the validity of properties is required only on conditional set, not on the whole set. The inspiration for its introduction came from the conditional expectation, which is the fundamental
notion of probability theory. Let us also remark that this operator generalizes the aggregation operator introduced earlier by Calvo et al. in~\cite[Definition 1]{CalvoKolesarovaKomornikovaMesiar2002} and it is the crucial ingredient in the definition of the generalized survival function.

\begin{definition}\label{def: cao}\rm 
(cf. \cite[Definition 3.1]{BoczekHalcinovaHutnikKaluszka2020})
A map $\aA[B]{\cdot}\colon [0,\infty)^{[n]}\to[0,\infty)$ is said to be a~\textit{conditional aggregation operator}  with respect to a set 
$B\in 2^{[n]}\setminus \{\emptyset$\} 
if it satisfies the following conditions:
\begin{itemize}
\item[i)] $\aA[B]{\x}\leq \aA[B]{\y}$ for any $\x,\y\in[0,\infty)^{[n]}$ such that $x_i\leq y_i$ for any $i\in B$;
\item[ii)] $\aA[B]{\bin_{B^c}}=0.$
\end{itemize}
\end{definition}

Let us compare the settings of the previous definition with the settings of the original definition, see~\cite[Definition 3.1]{BoczekHalcinovaHutnikKaluszka2020}). We consider the greatest $\sigma$-algebra as the domain of $\mu$ in comparison with the original arbitrary $\sigma$-algebra $\Sigma$. Then all vectors are measurable and this assumption may be omitted from the definition. The measurability of each vector is desired property mainly from the application point of view.
Because of the property $\aA[B]{\x}=\aA[B]{\x\bin_B}$ for any $\x\in[0,\infty)^{[n]}$ with fixed $B\in 2^{[n]}\setminus \{\emptyset\}$ the value $\aA[B]{\x}$ can be interpreted as \say{an aggregated value
of $\x$ on $B$}, see~\cite{BoczekHalcinovaHutnikKaluszka2020}.
In the following we list several examples of conditional aggregation operators we shall use in this paper. For further examples and some properties of conditional aggregation operators we recommend~\cite[Section~3]{BoczekHalcinovaHutnikKaluszka2020}.

\begin{example}\label{pr aggr}
Let $\x\in[0,\infty)^{[n]}$, $B\in 2^{[n]}\setminus\{\emptyset\}$ and $m\in\mathbf{M}$.
\begin{enumerate}[i)]
\item $\aAi[B]{\x}{m-\mathrm{ess}}=\mathrm{ess} \sup_m(\x\bin_{B})$, where $\mathrm{ess} \sup_{m}(\x)=\min\{\alpha\geq 0:\, \{\x>\alpha\}\in\mathcal{N}_{m}\}$.\footnote{A~set $N\in 2^{[n]}$ is said to be a~{\it null set} with respect to a monotone measure $m$ if $m(E\cup N)=m(E)$ for all $E\in 2^{[n]}.$ By $\cN$ we denote the family of null sets with respect to $m.$}
\item $\aA[B]{\x}=\mathrm{J}(\x\bin_{B},m)$, (the multiplication of vectors is meant by components) where $\mathrm{J}$ is an integral defined in~\cite[Definition 2.2]{BoczekHovanaHutnikKaluszka2020b}. 
Namely,\begin{itemize}
\item[a)]$\aAi[B]{\x}{\mathrm{Ch}_m}=\sum\limits_{i=1}^{n}x_{(i)}\left(m(\Ai{i}{\cap B})-m(\Ai{i+1}{\cap B})\right)$;
\item[b)]$\aAi[B]{\x}{\mathrm{Sh}_m}=\max\limits_{i\in[n]}\left\{ x_{(i)}\cdot m(\Ai{i}{\cap B} )\right\}$;
\item[c)]$\aAi[B]{\x}{\mathrm{Su}_m}=\max\limits_{i\in[n]}\left\{ \min\{x_{(i)}, m(\Ai{i}{\cap B})\} \right\}$.
\end{itemize}
\item $\aA[B]{\x}=
     \frac{\max_{i\in B}(x_i\cdot w_i)}{\max_{i\in B} z_i},$ where $\w\in[0,1]^{[n]}$ is a fixed weight vector, $\z\in(0,1]^{[n]}$ is fixed vector such that $\max_{i\in [n]}z_i=1$. We note, that for $\w=\z=\mathbf{1}_{[n]}$ we get $\aAi[B]{\x}{\mathrm{max}}=\max_{i\in B}x_i$.
\item
$\aAi[B]{\x}{p-\mathrm{mean}}=\left(\frac{1}{\counting(B)}\cdot\sum\limits_{i\in B}(x_i)^{p}\right)^{\frac{1}{p}}$ with $p\in(0,\infty)$. For $p=1$ we get the arithmetic mean.
\item $\aAi[B]{\x}{\mathrm{size}}=\max\limits_{D\in\collection}\size(\x\bin_{B})(D)$ with $\size$ being a size, see~\cite{BorzovaHalcinovaSupina2018, Halcinova2017, HalcinovaHutnikKiselakSupina2019}, is the outer essential supremum of $\x$ over $B$ with respect to a size $\size$ and a collection $\collection\subseteq2^{[n]}$. In particular, for the sum as a size, i.e., $\size_{\mathrm{sum}}(\x)(G)=\sum\limits_{i\in G}x_i$ for any $G\in2^{[n]}$ and for $\collection$ such that there is a set  $C\supseteq B, C\in\collection$ %containing $B$, resp. superset of $B$ 
we get
$\aAi[B]{\x}{\mathrm{sum}}=\sum\limits_{i\in B}x_i$.
\end{enumerate}
\end{example}

Observe that the empty set is not included in the Definition~\ref{def: cao}. The reason for that is the fact that the empty set does not provide any additional information for aggregation. However, in order to have the concept of the generalized survival function correctly introduced, it is necessary to add the assumption $\aA[\emptyset]{\cdot}=0$, see~\cite[Section 4]{BoczekHalcinovaHutnikKaluszka2020}. 
From now on, we shall consider only these conditional aggregation operators. 
Let us remark, that all mappings from Example~\ref{pr aggr} with the convention \enquote{$0/0=0$} satisfy this property. In the following we shall provide the definition of the generalized survival function, see~\cite[Definition 4.1.]{BoczekHalcinovaHutnikKaluszka2020}.
Let us consider a collection $\collectionGSF$,  ${{\{}}\emptyset{{\}}}\subseteq\collectionGSF\subseteq2^{[n]}$ and conditional aggregation operators on sets from $\collectionGSF$ with $ \aA[\emptyset]{\cdot}=0$.
The set of such aggregation operators we shall denote by $$\bA=\{\aA[E]{\cdot}: E\in \collectionGSF\}\footnote{Since $\aA[E]{\cdot}\colon[0,\infty)^{[n]}\to [0,\infty)$, $\cA$ is a~family of operators parametrized by a~set from~${\cE}$.}$$ and we shall call it a \textit{family of conditional aggregation operators} (FCA for short). For example, $\cA^{\mathrm{sum}}=\{\sA^{\mathrm{sum}}(\cdot|E):E\in2^{[n]}\}$,
$\cA^{\mathrm{max}}=\{\sA^{\mathrm{max}}(\cdot|E):E\in\{\emptyset, \{1\}, \{2\}, \dots, \{n\}\}\}$, 
$\widehat{\cA}^{\mathrm{max}}=\{\sA^{\mathrm{max}}(\cdot|E):E\in\{\emptyset\}\}$
or $\cA=\{\sA(\cdot|E):E\in2^{[n]}\}$, $n\geq 2$, where
$$\sA(\x|E)=\begin{cases}
\sA^{\mathrm{max}}(\x|E), & E\in\{\{1\}, \{2\}, \dots, \{n\}\},\\
0,&E=\emptyset,\\
\sA^{\mathrm{sum}}(\x|E), &\text{otherwise}
\end{cases}$$
for any $\x\in[0,\infty)^{[n]}$.

\begin{definition}\label{def: gsf}\rm(cf.~\cite[Definition 4.1.]{BoczekHalcinovaHutnikKaluszka2020})
Let  $\x\in[0,\infty)^{[n]}$, $\mu\in\mathbf{M}$. The \textit{generalized survival function} with respect to $\bA$ is defined as
\begin{eqnarray*}%\label{predpis gsf}
\gsf{}{\x}{\alpha}=\min\left\{\mu(E^c): \aA[E]{\x}\leq\alpha,\, E\in\collectionGSF\right\}\end{eqnarray*} for any $\alpha\in[0,\infty)$. 
\end{definition}

The presented definition is correct. Really, for any $E\in\collectionGSF$ it holds that $E^c\in2^{[n]}$ is a~measurable set. Moreover, the set $\{E\in\collectionGSF: \aA[E]{\x}\leq\alpha\}$ is nonempty for all $\alpha\in[0,\infty)$, because $\aA[\emptyset]{\cdot}=0$ by convention and $\emptyset\in\collectionGSF$. Immediately it is seen, that for $\collectionGSF=2^{[n]}$ and $\cA^{\mathrm{max}}$ we get the standard survival function, compare with~\eqref{prepis2}.
When it will be necessary we shall emphasize the collection $\collectionGSF$ in the notation of generalized survival function, i.e. we shall use $\cA^{\collectionGSF}$.

On several places in this paper we shall work with the FCA that is \textit{nondecreasing} w.r.t sets, i.e. the map $E\mapsto\aA[E]{\cdot}$ will be nondecreasing. Many FCA satisfy this property, e.g. $\cA^{m-\mathrm{ess}}=\{\sA^{m-\mathrm{ess}}(\cdot|E):E\in\collectionGSF\}$, $\cA^{\mathrm{Ch}_m}=\{\sA^{\mathrm{Ch}_m}(\cdot|E):E\in\collectionGSF\}$, $\cA^{\mathrm{Su}_m}=\{\sA^{\mathrm{Su}_m}(\cdot|E):E\in\collectionGSF\}$, $\cA^{\mathrm{Sh}_m}=\{\sA^{\mathrm{Sh}_m}(\cdot|E):E\in\collectionGSF\}$, $\cA^{\mathrm{max}}=\{\sA^{\mathrm{max}}(\cdot|E):E\in\collectionGSF\}$, see Example~\ref{pr aggr} i), ii), iii).

\section{Equality and inequalities of the generalized and standard survival function}

In this section we shall treat Problem 1. We provide sufficient and necessary conditions on $\x$, $\mu$ and $\cA$ under which the generalized survival function and survival function coincide. 
%  The motivation for studying this problem is the observation that survival functions can coincide for specific monotone measures, vectors and FCA, see the Introduction. 
The important knowledge we use is the standard survival function formula.
In what follows we shall work with the expression of the survival function on a finite set in the form 
% \begin{equation}\label{LM_form}\lm=\sum_{i=1}^n\mu\left(\Ai{i}\right)\cdot\mathbf{1}_{[x_{(i-1)},x_{(i)})}(\alpha)\text{,}\end{equation}
\begin{equation}\label{LM_form}\lm=\sum_{i=0}^{n-1}\mu\left(\Ai{i+1}\right)\cdot\mathbf{1}_{[x_{(i)},x_{(i+1)})}(\alpha)\text{}\end{equation}
with the permutation $(\cdot)$ such that $0=x_{(0)}\leq x_{(1)}\leq x_{(2)}\leq\dots \leq x_{(n)}$ and
$\Ai{i}=\{(i),\dots,(n)\}$ for $i\in[n]$. %{\lh{Preco tu pouzivame kcka? a nie icka?}}, 
 % see e.g. the Choquet integral formula on discrete set~\cite{ChenMesiarLiStupnanova2017}. 
% \lh{??? POUZIT TVAR:  
% \begin{equation}\label{LM_form*}\lm=\sum_{i=0}^{n-1}\mu\left(\Ai{i+1}\right)\cdot\mathbf{1}_{[x_{(i)},x_{(i+1)})}(\alpha)\text{,}\end{equation}
% with the permutation $(\cdot)$ and
% $\Ai{i+1}=\{(i+1),\dots,(n)\}$ for $i\in[n]\cup\{0\}$ with the convention $\Ai{n+1}=\emptyset$.}
However, one can easily see that some summands in the formula~\eqref{LM_form} can be redundant. For example, for vectors with the property $x_{(i)}=x_{(i+1)}$ for some $i\in[n-1]\cup\{0\}$ we have $\mu\left(\Ai{i+1}\right)\cdot\mathbf{1}_{[x_{(i)},x_{(i+1)})}(\alpha)=0$ for any $\alpha\in[0,\infty)$, i.e., this summand does not change the values of survival function and can be omitted.
\bigskip

Let us consider an arbitrary (fixed)  input vector $\x$ together with a
% (fixed)
permutation $(\cdot)$ such that $0=x_{(0)}\leq x_{(1)}\leq x_{(2)}\leq\dots \leq x_{(n)}$. %Henceforward, we shall not explicitly mention this permutation in the assumptions of presented results. 
Let us denote \begin{align}\label{pi}\Psi_{\x}:=\{i\in[n-1]\cup\{0\}: x_{(i)}<x_{(i+1)}\}\cup\{n\}.\end{align}
For example, for the input vector $\x=(3,2,3,1)$ and the permutation $(\cdot)$ such that $x_{(0)}=0$, $x_{(1)}=1$, $x_{(2)}=2$, $x_{(3)}=3$, $x_{(4)}=3$, we get 
$\Psi_{\x}=\{0,1,2,4\}$. The following proposition includes the very basic properties of system $\Psi_{\x}$ needed for further results.
\begin{proposition}\label{vlastnostipi}
Let $\x\in[0,\infty)^{[n]}$.
% Let $\Psi_{\x}$ be a system given by~(\ref{pi}) with $(\cdot)$ being a permutation of $\x$ such that $0=x_{(0)}\leq x_{(1)}\leq x_{(2)}\leq\dots \leq x_{(n)}$. 
\begin{enumerate}[i)]
\item $\Psi_{\x}$ is independent on permutation $(\cdot)$ of $\x$,  i.e., $\Psi_{\x}$ contains the same values for any permutation $(\cdot)$ of $\x$ such that $0=x_{(0)}\leq x_{(1)}\leq x_{(2)}\leq\dots \leq x_{(n)}$.
\item For any $i\in[n]$ there exists $k_i\in\Psi_{\x}\setminus\{0\}$ such that $x_i=x_{(k_i)}$, i.e. $\{x_{(k_i)}:k_i\in\Psi_{\x}\setminus\{0\}\}$ contains all different values of $\x$.
\item $x_{(\min\Psi_{\x})}=0$.
\item $\left\{[x_{(k)},x_{(k+1)}): k\in\Psi_{\x}\right\}$ is a decomposition of interval $[0,\infty)$ into nonempty pairwise disjoint~sets.
\end{enumerate}
\end{proposition}
\proof\begin{enumerate}[i)]
\item 
% Let us consider two different nondecreasing permutations of $\x$ (if they exist) $(\cdot)_1$ and $(\cdot)_2$. Let us set $\pi_1:=\{i\in[n-1]\cup\{0\}: x_{(i)_1}<x_{(i+1)_1}\}\cup\{n\}$ and $\pi_2:=\{i\in[n-1]\cup\{0\}: x_{(i)_2}<x_{(i+1)_2}\}\cup\{n\}$. We shall show that $\pi_1=\pi_2$.

% At first, let us take $i=0$. If $0=x_{(0)_1}<x_{(1)_1}$, then there is not $i\in[n]$ such that $x_i$=0 and it follows that $x_{(0)_2}<x_{(1)_2}$. This means that $0\in\pi_1$ and $0\in\pi_2$. On the other hand, if $x_{(0)_1}=x_{(1)_1}$, then there exists $i\in[n]$ such that $x_i=0$. Then it holds that $x_{(0)_2}=x_{(1)_2}$, what means that neither $0\notin\pi_1$ nor $0\notin\pi_2$.

% Let us suppose that $\{i\in[n-1]: x_{(i)_1}<x_{(i+1)_1}\}\neq\{i\in[n-1]: x_{(i)_2}<x_{(i+1)_2}\}$. Without loss of generality, we may suppose that there exists $k_i\in[n-1]$ such that $x_{(k_i)_1}<x_{(k_i+1)_1}$ and $x_{(k_i)_2}=x_{(k_i+1)_2}$. Further, from the nondecreasingness of permutations $(\cdot)_1$, resp. $(\cdot)_ 2$ it holds for any $i\in[n]$ that $\{j_i\in[n]:x_i=x_{(j_i)_1}\}=\{j_i\in[n]:x_i=x_{(j_i)_2}\}$. For the completeness of the proof it is enough to take $i\in[n]$ such that $x_i=x_{(k_i)_2}$. Then $k_i, k_i+1 \in\{j_i\in[n]:x_i=x_{(j_i)_2}\}=\{j_i\in[n]:x_i=x_{(j_i)_1}\}$, what is a contradiction because $x_{(k_i)_1}<x_{(k_i+1)_1}$.
Let us consider two different permutations of $\x$ (if they exist) $(\cdot)_1$ and $(\cdot)_2$ with the required property. Let us denote
\begin{align*}\Psi_\x&:=\{i\in[n-1]\cup\{0\}: x_{(i)_1}<x_{(i+1)_1}\}\cup\{n\},\\\Phi_{\x}&:=\{i\in[n-1]\cup\{0\}: x_{(i)_2}<x_{(i+1)_2}\}\cup\{n\}.
\end{align*}
We show that $\Psi_{\x}=\Phi_{\x}$. Indeed, $n\in\Psi_{\x},n\in\Phi_{\x}$. If $i\in\Psi_{\x}\setminus\{n\}$, then $x_{(i)_1}<x_{(i+1)_1}$. Because of nondecreasing rearangement of $\x$ with respect to $(\cdot)_1$, $(\cdot)_2$ we get $x_{(i)_2}<x_{(i+1)_2}$, therefore  $i\in\Phi_{\x}$ and $\Psi_{\x}\subseteq\Phi_{\x}$. By analogy it holds $\Phi_{\x}\subseteq\Psi_{\x}$.
\item
% We shall show that for any $i\in[n]$ there exists $k_i\in\Psi_{\x}\setminus\{0\}$ such that $x_i=x_{(k_i)}$. 
Since any $i\in[n]$ can be represented via permutation as $i={(j_i)}$, 
% where
$j_i\in[n]$, let us set $$k_i=\max\{j_i\in[n]:\, x_i=x_{(j_i)}\}.$$
As for any $k_i<n$ it holds that $x_{(k_i)}<x_{(k_{i}+1)}$, then $k_i\in\Psi_{\x}\setminus\{0\}$. Moreover, $k_i=n\in\Psi_{\x}$ because of the definition of $\Psi_{\x}$, see~\eqref{pi}.

\item  It follows immediately from the fact that $\min\Psi_{\x}=\max\{i\in[n]\cup\{0\}:x_{(i)}=x_{(0)}=0\}$.
% {\jb{If $x_{(1)}\neq 0$, then $x_{(0)}<x_{(1)}$ and hence $\min\Psi_{\x}=0$. Consequently, $x_{(\min\Psi_{\x})}=0$. On the other hand, if $x_{(1)}=0$, then $\min\Psi_{\x}=\max\{j_i\in[n]:0=x_{(j_i)}\}$. It follows that $x_{(\min\Psi_{\x})}=0$.}}

\item It follows from part ii), iii) and from definition of system $\Psi_{\x}$, since  $x_{(k)}<x_{(k+1)}$ for any $k\in\Psi_{\x}$ and $x_{(k_1)}\neq x_{(k_2)}$ for any $k_1,k_2\in\Psi_{\x}$.
% {\color{orange}
% alebo: From definition of system $\Psi_{\x}$ we have that $[0,\infty)=\left(\bigcup_{k\in\Psi_{\x}}[x_{(k)},x_{(k+1)})\right)\cup\left(\bigcup_{k\in[n]\setminus\Psi_{\x}}[x_{(k)},x_{(k+1)})\right)=\left(\bigcup_{k\in\Psi_{\x}}[x_{(k)},x_{(k+1)})\right)\cup\bigcup_{k\in[n]\setminus\Psi_{\x}}\emptyset=\bigcup_{k\in\Psi_{\x}}[x_{(k)},x_{(k+1)})$.
% }
% \item {\bs{It is immediate.
% % The inequality $\aAi[E]{\x}{\mathrm{max}}<x_{(k+1)}$, $k\in\Psi_{\x}$, means $E\subseteq\{x_{(1)},\dots,x_{(k)}\}=\Ai{k+1}^c$.
% }}
\qed
\end{enumerate}

Since we have shown that the system $\Psi_{\x}$ is independent of the chosen permutation, henceforward we shall not explicitly mention the permutation in assumptions of presented results. The following proposition states that the formula~(\ref{LM_form}) can be rewritten by the system $\Psi_{\x}$ in the simpler form, see part ${\rm{i)}}$.
% Moreover, in the second part of the proposition we provide a formula for the generalized survival function w.r.t. $\cA^{\mathrm{max}}$ with the collection dependent on $\Psi_{\x}$.
Moreover, in the second part of the proposition we show that {for a fixed $\x\in[0,\infty)^{[n]}$ it is $\gsf{\mathrm{max}, {{\collectionGSF}}}{\x}{\alpha}= \lm$}  with smaller collection $\collectionGSF$ than the whole powerset $2^{[n]}$ (compare with the known result~\cite[Example 4.2]{BoczekHalcinovaHutnikKaluszka2020} or see \eqref{prepis2}). The collection $\collectionGSF$ depends on $\x$ (equivalently on $\Psi_{\x}$).
% With regard to above written facts, the following result is true. The result says that for the family of conditional aggregation operators $\cA^{\mathrm{max}}$ the (generalized and standard) survival functions equal also for smaller system of sets than the powerset $2^{[n]}$ (compare with the known result ~\cite[Example 4.2]{BoczekHalcinovaHutnikKaluszka2020} or see \eqref{prepis2}).

\begin{proposition}\label{LM=SLM_max}
Let $\x\in[0,\infty)^{[n]}$, $\mu\in\mathbf{M}$.
\begin{enumerate}[i)]
    \item Then \begin{equation}\label{PiLM_form}\lm=\sum_{k\in\Psi_{\x}}\mu\left(\Ai{k+1}\right)\cdot\mathbf{1}_{[x_{(k)},x_{(k+1)})}(\alpha)\end{equation}
    for any $\alpha\in[0,\infty)$ with the convention $x_{(n+1)}=\infty$. %\label{PiLM_form}
    \item If $\collectionGSF\supseteq\{\Ai{k+1}^c: k\in\Psi_{\x}\}$, then $$\gsf{\mathrm{max}}{\x}{\alpha}= \lm.$$ 
    %\label{LM=SLM_max}
\end{enumerate}
\end{proposition}
\proof
\begin{enumerate}[i)]
    \item For $k\in[n-1]\cup\{0\},\, k\notin\Psi_{\x}$, we have $x_{(k)}=x_{(k+1)}$. 
    This leads to the fact that  $\mu\left(\Ai{k+1}\right)\cdot\mathbf{1}_{[x_{(k)},x_{(k+1)})}(\alpha)=0$ for any $\alpha\in[0,\infty)$. {Using the Proposition~\ref{vlastnostipi} iv) we have the required assertion.}
    \item According to Proposition~\ref{vlastnostipi} (iv) let us divide interval $[0,\infty)$ into disjoint sets $$[0,\infty)=\bigcup_{k\in\Psi_{\x}}[x_{(k)},x_{(k+1)}).$$
    Let us consider an arbitrary (fixed) $k\in\Psi_{\x}$.
    Then from the fact, that $\Ai{k+1}^c\in\collectionGSF$ and $\aAi[E_{(k+1)}^c]{\x}{\mathrm{max}}=x_{(k)}$
    we have {$E_{(k+1)}^c\in\{E:\aA[E]{\x}\leq\alpha\}$ for any $\alpha\in[x_{(k)},x_{(k+1)})$. } 
    %$\mu(E_{(k+1)}^c)\in\{\mu(E^c): \aAi[E]{\x}{\mathrm{max}}\leq\alpha,E\in\collectionGSF\}$. 
    Therefore we get $$\gsf{\mathrm{max}, {{\collectionGSF}}}{\x}{\alpha}=\min\{\mu(E^c):\aAi[E]{\x}{\mathrm{max}}\leq\alpha,E\in\collectionGSF\}\leq\mu(E_{(k+1)})=\lm$$
    for any $\alpha\in[x_{(k)},x_{(k+1)})$, where the last equality follows from part~$\text{i)}$.
    On the other hand, {as $\collectionGSF\subseteq2^{[n]}$ from properties of minimum we have} $$\gsf{\mathrm{max}, {{\collectionGSF}}}{\x}{\alpha}\geq\gsf{\mathrm{max}, {{2^{[n]}}}}{\x}{\alpha}=\mu(E_{(k+1)})=\lm$$ for any $\alpha\in[x_{(k)},x_{(k+1)})$. %, where first inequality holds, because $\collectionGSF\subseteq2^{[n]}$.
    To sum it up, $\gsf{\mathrm{max}, {{\collectionGSF}}}{\x}{\alpha}=\mu(E_{(k+1)})=\lm$ for any $\alpha\in[x_{(k)},x_{(k+1)})$.\qed
\end{enumerate}

\begin{remark}
Let us remark that in the whole paper we suppose $\mu$ is defined on $2^{[n]}$. However, in fact, it is enough to have a smaller $\mathtt{Dom}{(\mu)}$. For example,  in part ii) of the previous proposition it is enough to consider the domain of $\mu$ being $\{E^c: E\in\mathcal{E}\}$.
\end{remark}

Let us note that in~(\ref{PiLM_form}) the last summand is always equal to $0$ because $\mu\left(\Ai{n+1}\right)=\mu(\emptyset)=0$. However, it is useful to consider the form of survival function in~(\ref{PiLM_form}) with sum over the whole set $\Psi_{\x}$ not $\Psi_{\x}\setminus\{n\}$  because of some technical details in presented proofs in this paper.

%\begin{remark}\label{remark}
%Let us remark that the set $\{x_{(k_{i})}:\, k_{i}\in\Psi_{\x}\}$ covers all values of input vector $\x$, i.e., for any $i\in[n]$ there exists $k_i\in\Psi_{\x}$ such that $x_i=x_{(k_i)}$. Indeed, any $i\in[n]$  can be expressed via permutation as $i=(j_i)$. Then it is enough to take 
%$$k_i=\max\{j_i\in[n]:\, x_i=x_{(j_i)}\}.$$
%It is easy to see that $k_i\in\Psi_{\x}$, since $x_{(k_i)}<x_{(k_{i}+1)}$ for $k_i< n$ and if $k_i=n$, then $k_i\in\Psi_{\x}$ from the definition of $\Psi_{\x}$.
%Moreover, this set always contains the zero value (because of the definition of $\Psi_{\x}$), although $0$ need not be one of components of the input vector.
%More precisely, if there is no component of input vector with zero value, then $0\in\Psi_{\x}$ and $x_{(0)}=0$. Otherwise, $x_{(\min \Psi_{\x})}=0$.
%\end{remark}

% It will turn out later that this generalization will help us to formulate crucial results.

% \begin{align*}
% \lm&=\mu([n]\setminus \{i\in [n]: x_i\leq \alpha\})\\&=\min\big\{\mu (E^c):(\forall i\in E)\,\, x_i\leq \alpha,\, E\in \collectionGSF\supseteq\{\Ai{k+1}^c: k\in\Psi_{\x}\}\big\}\nonumber\\&=\min\big\{\mu (E^c): \aAi[\x]{E}{\mathrm{max}}\leq \alpha,\, E\in \collectionGSF\supseteq\{\Ai{k+1}^c: k\in\Psi_{\x}\}\big\},
% \end{align*}

\begin{example}\label{LM_collection} Let us take $\cA^{\mathrm{max}}=\{\aAi[E]{\cdot}{\mathrm{max}}: E\in\collectionGSF\}$ and normalized monotone measure $\mu$ on $2^{[3]}$ given in the following table:
\begin{table}[H]
\renewcommand*{\arraystretch}{1.2}
\begin{center}
\begin{tabular}{|c|c|c|c|c|c|c|c|c|}
\hline
$E$ & $\emptyset$ & $\{1\}$ & $\{2\}$ & $\{3\}$ & $\{1,2\}$ & $\{1,3\}$ & $\{2,3\}$ & $\{1,2,3\}$ \\ \hline
$\mu(E)$ & $0$ & $0$ & $0.5$ & $0$ & $0.5$ & $0.5$ & $0.5$ & $1$\\ \hline
\end{tabular}
%\caption{The values of normalized monotone measure from Example~\ref{LM_collection}}
\label{table0}
\end{center}
\vspace*{-12pt}
\end{table}
% $\mu(\emptyset)=\mu(\{1\})=\mu(\{3\})=0$, $\mu(\{2\})=\mu(\{1,2\})=\mu(\{2,3\})=\mu(\{1,3\})=0{.}5$. 
\noindent Further, let us take the input vector $\x=(1,2,1)$ with the permutation \mbox{$(1)\!=\!1$, $(2)\!=\!3$, $(3)\!=\!2$}.

Then $\Psi_{\x}=\{0,2,3\}$ and the collection guarantying the equality between survival function and generalized survival function (of input $\x$) is according to Proposition~\ref{LM=SLM_max} $\text{ii)}$ e.g.
\begin{align*}
\collectionGSF&=\{\Ai{k+1}^c:\, k\in\Psi_{\x}\}=\{\Ai{1}^c, \Ai{3}^c, \Ai{4}^c\}=\{\emptyset, \{(1), (2)\}, \{(1), (2), (3)\}\}\\&=\{\emptyset, \{1,3\}, \{1,2,3\}\}.\end{align*}
% From Proposition~\ref{LM=SLM_max}, part $\text{ii)}$ we have
% $$\lm=\gsf{\mathrm{max}, {{\collectionGSF}}}{\x}{\alpha}=\min\{\mu(E^c):\, \aAi[E]{\x}{\mathrm{max}}\leq\alpha, \, E\in\{\Ai{k+1}^c: k\in\Psi_{\x}\}\},$$
% what is really true, since 
Indeed, 
$$\gsf{\mathrm{max}}{\x}{\alpha}=
% 0\cdot\bin_{[2,\infty)}+ 0.5\cdot\bin_{[1,2)}+1\cdot\bin_{[0,1)}=
1\cdot\bin_{[0,1)}(\alpha)+0.5\cdot\bin_{[1,2)}(\alpha)=\lm.$$
\end{example}

\begin{figure}
  \begin{center}
   \includegraphics[scale=1]{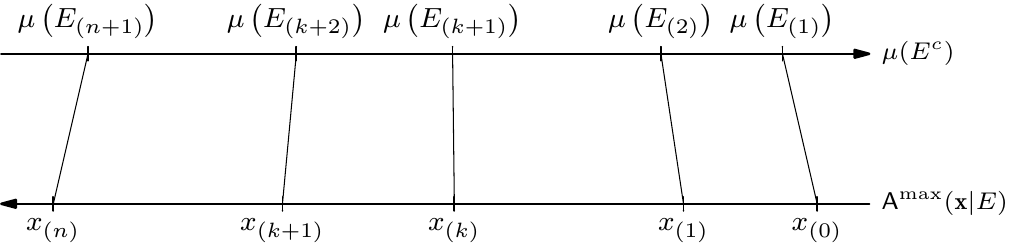}
   \caption{The survival function visualization $\mu_{\cA^{\mathrm{max}}}$ with $\cA=\{\aAi[E]{\cdot}{\mathrm{max}}: E\in\{\Ai{k+1}^c: k\in\Psi_{\x}\}\}$}
   \label{LM_visualization}
    \end{center}
\end{figure}

From the previous result it follows that the standard survival function can be represented by the formula 
\begin{align}\label{formula}
\lm&=\min\big\{\mu (E^c): \aAi[E]{\x}{\mathrm{max}}\leq \alpha,\, E\in \{\Ai{k+1}^c: k\in\Psi_{\x}\}\big\}
\end{align}
 with the system $\Psi_{\x}$
%  permutation $(\cdot)$ 
 given by the input vector $\x$. This formula can be visualized by Figure~\ref{LM_visualization}. Let us remark that since $\aAi[\Ai{k+1}^c]{\x}{\mathrm{max}}=x_{(k)}$ on the upper line we measure sets $(\Ai{k+1}^c)^c$. The calculation of (generalized) survival function is processed as we have described in the Introduction. 
%  At the end of this section let us remark that the essence of this paper is
Let us remark that the essence of the following results is 
the pointwise comparison of the generalized survival function with the standard survival function 
% %  given by~\eqref{PiLM_form} 
having in mind the representation~\eqref{formula} together with its visualization, see Figure~\ref{LM_visualization}.

It is obvious that the equality of survival functions (standard and generalized) means that they have to achieve the same values, i.e., $\mu\left(\Ai{k+1}\right)$, $k\in\Psi_{\x}$, on the same corresponding intervals $[x_{(k)},x_{(k+1)})$, $k\in\Psi_{\x}$.
%It is clear that if survival functions (standard and generalized) equal for the given input $\x$, then the generalized survival function has to achieve the same values as the survival function, i.e., $\mu\left(\Ai{k+1}\right)$, $k\in\Psi_{\x}$ and these values have to be achieved on the corresponding intervals $[x_{(k)},x_{(k+1)})$, $k\in\Psi_{\x}$. 
Having in mind the formula~\eqref{PiLM_form}, the survival function representation given by~\eqref{formula} and the visualization, see Figure~\ref{LM_visualization}, we can formulate the following sufficient conditions.
% for the equality between survival functions. 
While 
% the condition 
(C1) ensures that the generalized survival function will be able to achieve the same values as the survival function, 
% the condition 
(C2) guarantees it. Let $\cA$ be FCA.
%  Condition (C3) is a slighter version of (C1), which will be also useful.  Let $\x\in[0,+\infty)^{[n]}$ with the permutation $(\cdot)$.
\begin{enumerate}[(C1)]
    \item For any $k\in\Psi_{\x}$ there exists $G_k\in\collectionGSF$ such that $$
    \aA[G_k]{\x}=x_{(k)}\quad\text{and}\quad \mu(G_k^c)=\mu(\Ai{k+1}).$$ 
    \item For any $k\in\Psi_{\x}$ and for any  $E\in \collectionGSF$ it holds:
    $$\aA[E]{\x}<x_{(k+1)}\Rightarrow\mu(E^c)\geq \mu(\Ai{k+1}).$$
  \end{enumerate}
 
  \begin{figure}
  \begin{center}
   \includegraphics[scale=1.1]{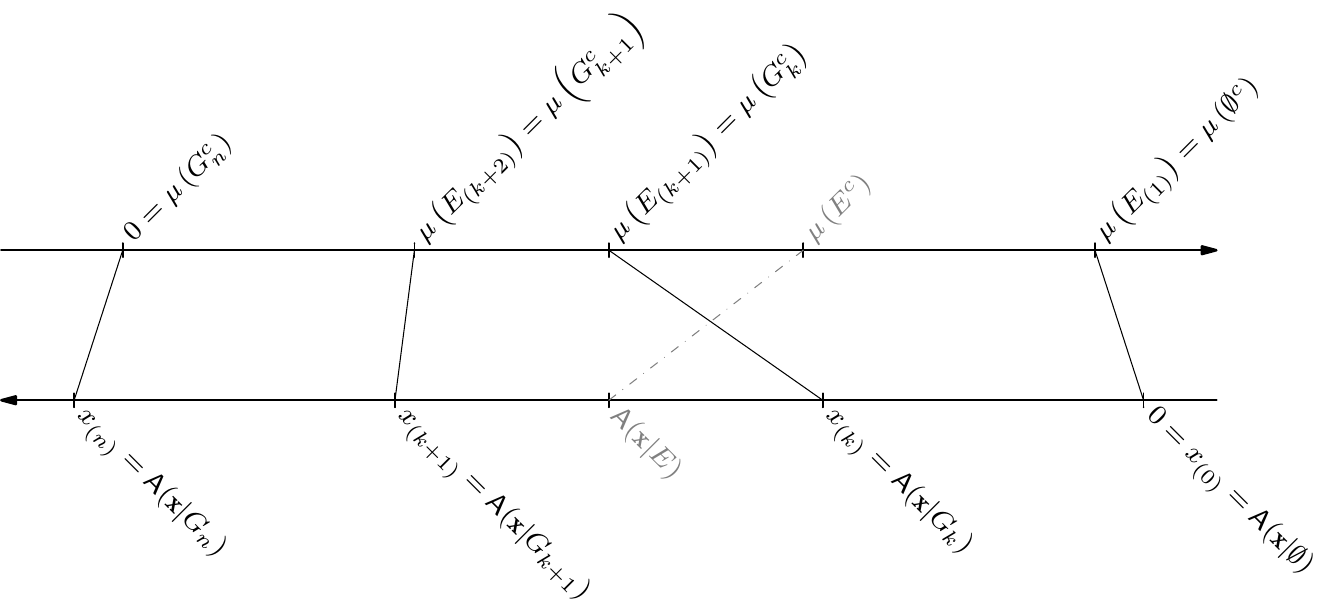}
   \caption{The visualization of conditions $\mathrm{(C1)}$ and $\mathrm{(C2)}$}
   \label{obr_dokaz}
    \end{center}
\end{figure}
 
\noindent The visualization of conditions $\mathrm{(C1)}$, $\mathrm{(C2)}$ via two parallel lines is drawn in Figure~\ref{obr_dokaz}. Let us remark that for $k=n$ (C2) holds trivially. Also, for $k=\min\Psi_{\x}$ (C1) holds trivially with $G_{\min\Psi_{\x}}=\emptyset$.

\begin{remark}\label{poznamka_max} In accordance with the above written, it can be easily seen that for $\cA^{\mathrm{max}}$ with $\collectionGSF\supseteq\{\Ai{k+1}^c: k\in\Psi_{\x}\}$
it holds $G_k=\Ai{k+1}^c$ for any $k\in\Psi_{\x}$ regardless of the choice of $\mu$ in 
% the condition
$\mathrm{(C1)}$. Of course, for specific classes of monotone measures $\mu$ also other sets $G_k$ can  satisfy $\rm{(C1)}$. Similarly, the validity of 
% the condition
$\rm{(C2)}$ is clear. 
Indeed, if $\aAi[E]{\x}{\mathrm{max}}<x_{(k+1)}$, then  we have  $E\subseteq \Ai{k+1}^c$, i.e., $E^c\supseteq E_{(k+1)}$. From the monotonicity of $\mu$ we have $\mu(E^c)\geq\mu(E_{(k+1)})$ for any $E\in\collectionGSF$.
\end{remark}

Conditions {$\mathrm{(C1)}$} and {$\mathrm{(C2)}$} guarantee inequalities between survival functions. Thus the equality of survival function is a consequence.
% Having in mind conditions (C1), (C2)
% % and (C3)
% {\lh{and the proof of Corollary~\ref{LM_SLM_coincide}}} we , can naturally formulate results leading to inequalities between (generalized and standard) survival functions.

\begin{proposition}\label{LM<SLM}
Let $\x\in[0,\infty)^{[n]}$, $\mu\in\bM$, and let $\cA$ be FCA.
\begin{enumerate}[i)]
\item If $\mathrm{(C1)}$ holds, then  $\gsf{}{\x}{\alpha}\leq\lm$
for any $\alpha\in[0,\infty)$. %\label{LM>SLM_0}
\item $\mathrm{(C2)}$ holds if and only if $\lm\leq \gsf{}{\x}{\alpha}$
for any $\alpha\in[0,\infty)$.
\end{enumerate}
\end{proposition}
\proof According to Proposition~\ref{vlastnostipi} (iv) let us divide interval $[0,\infty)$ into disjoint sets $$[0,\infty)=\bigcup_{k\in\Psi_{\x}}[x_{(k)},x_{(k+1)}).$$ 
%with the convention $x_{(n+1)}=\infty$. 
% Obviously, from the definition of system $\Psi_{\x}$ intervals $[x_{(k)},x_{(k+1)})$, $k\in\Psi_{\x}$ are nonempty.
Let us consider an arbitrary (fixed) $k\in\Psi_{\x}$. 

Let us prove part i). According to
% condition 
{$\mathrm{(C1)}$} there exists the set $G_k\in\collectionGSF$ such that $\aA[G_k]{\x}=x_{(k)}$ and $\mu(G_k^c)=\mu(E_{(k+1)})$. From the fact that $\mu(E_{(k+1)})=\mu(G_k^c)\in\left\{\mu(E^c):\, \aA[E]{\x}\leq x_{(k)}\right\}$
 and since $\gsf{}{\x}{\alpha}$ is nonincreasing (see~\cite[Proposition 4.3 (a)]{BoczekHalcinovaHutnikKaluszka2020}) we have
$$\gsf{}{\x}{\alpha}\leq\gsf{}{\x}{x_{(k)}}\leq\mu(E_{(k+1)})=\lm$$
% $$\mu(E_{(k+1)})\geq\gsf{}{\x}{x_{(k)}}\geq\gsf{}{\x}{\alpha}$$
for any $\alpha\in[x_{(k)},x_{(k+1)})$, where the last equality follows from~(\ref{PiLM_form}).

Let us prove part ii). From 
% condition 
{$\mathrm{(C2)}$} it follows that for any $E\in\collectionGSF$ if $\aA[E]{\x}< x_{(k+1)}$, then $\mu(E^c)\geq\mu(E_{(k+1)})$. Therefore, $$\gsf{}{\x}{\alpha}\geq \mu(E_{(k+1)})=\lm$$ for any $\alpha\in[x_{(k)},x_{(k+1)})$, where the last equality follows from~(\ref{PiLM_form}).
It is enough to prove the implication  $\Leftarrow$.
% Let  $\lm\leq\gsf{}{\x}{\alpha}$ for any $\alpha\in[0,\infty)$, then {\lh{the inequality}} also holds for {\lh{any partial interval,}}  $\alpha\in[x_{(k)},x_{(k+1)})$ with $k\in\Psi_{\x}$. 
% From this fact and from~\eqref{PiLM_form} it holds
Since 
% $\lm\leq\gsf{}{\x}{\alpha}$
$$\gsf{}{\x}{\alpha}=\min\left\{\mu(E^c):\aA[E]{\x}\leq\alpha<x_{(k+1)}, E\in\collectionGSF\right\}\geq\mu(\Ai{k+1})=\lm$$
for any $\alpha\in[x_{(k)},x_{(k+1)})$, then for any $E\in\collectionGSF$ it has to hold: if $\aA[E]{\x}<x_{(k+1)}$, then $\mu(E^c)\geq\mu(\Ai{k+1})$. %Otherwise we get the contradiction.
\qed
\bigskip

\begin{corollary}\label{LM_SLM_coincide}
Let $\x\in[0,\infty)^{[n]}$, $\mu\in\bM$, and let $\cA$ be FCA. If 
% conditions
{$\mathrm{(C1)}$} and {$\mathrm{(C2)}$} are satisfied, then $\gsf{}{\x}{\alpha}= \lm$
for any $\alpha\in[0,\infty)$.
\end{corollary}

% \begin{remark}
% % Let us remark that the set $G_0$ in the condition $\mathrm{(C1)}$ always exists and can be identified with $\emptyset$. Indeed, $\aA[\emptyset]{\x}=0$ and $\mu([n])=\mu(\Ai{1})$. 
% It is enough to require the existence of sets $G_k$, $k\in\Psi_{\x}$ in the previous proposition. Any other sets would be additional. TBA
% \end{remark}

The application of the previous result is illustrated in the following example. The second example proves that $\mathrm{(C1)}$ and $\mathrm{(C2)}$ are only sufficient and not necessary.
% for the identity of generalized and survival functions in general.

\begin{figure}
\begin{center}
\begin{tabular}{m{11cm} m{4.3cm}}
\includegraphics[scale=0.9]{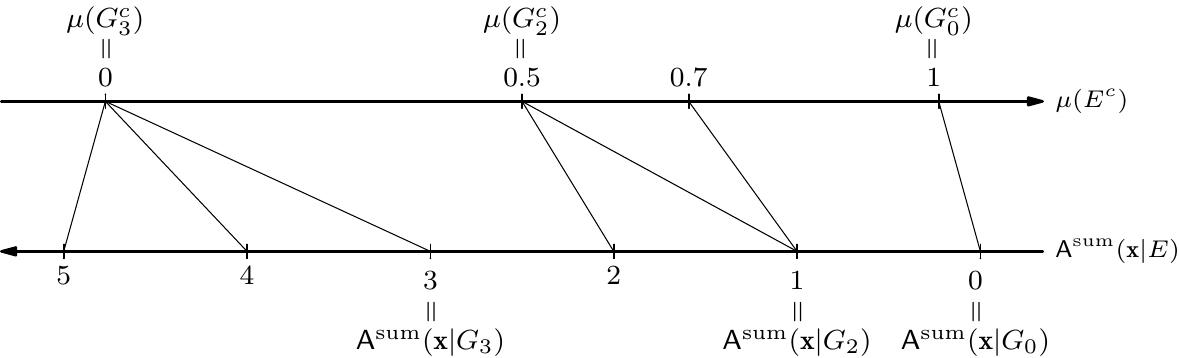} &
{\vspace*{-0.25cm}\includegraphics[scale=0.9]{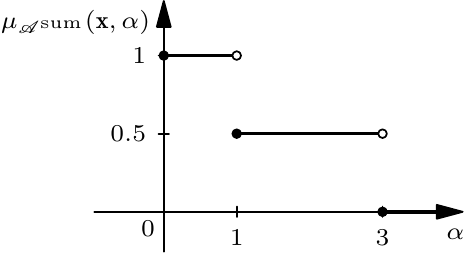}}
\end{tabular}
\caption{Generalized survival function and visualization from Example~\ref{SC_example}}  \label{LM_eq_SLM_picture}
   \end{center}
\end{figure}

\begin{example}\label{SC_example}
Let us consider $\cA^{\mathrm{sum}}=\{\aAi[E]{\cdot}{\mathrm{sum}}: E\in2^{[3]}\}$, and normalized monotone measure $\mu$ on $2^{[3]}$  with the following values:
% Let us consider $n=3$, $\collectionGSF=2^{[3]}$, the input vector $\x=(1,3,1)$ with the permutation $(1)=1$, $(2)=3$, $(3)=2$, $\sA^{\mathrm{sum}}$ and the normalized monotone measure $\mu$ on $2^{[3]}$ with the following values:
\begin{table}[H]
\renewcommand*{\arraystretch}{1.2}
\begin{center}
\begin{tabular}{|c|c|c|c|c|c|c|c|c|}
\hline
$E$ & $\emptyset$ & $\{1\}$ & $\{2\}$ & $\{3\}$ & $\{1,2\}$ & $\{1,3\}$ & $\{2,3\}$ & $\{1,2,3\}$ \\ \hline
$\mu(E)$ & $0$ & $0$ & $0.5$ & $0$ & $0.5$ & $0$ & $0.7$ & $1$\\ \hline
 $\aAi[E]{\x}{\mathrm{sum}}$ & $0$  & $1$ & $3$ & $1$ & $4$ & $2$ & $4$ & $5$\\\hline
\end{tabular}
%\caption{The values of conditional aggregation operator $\sA^{\mathrm{sum}}$ from Example~\ref{SC_example}}
%\label{table}
\end{center}
\vspace*{-12pt}
\end{table}
\noindent Further, let us take the input vector $\x=(1,3,1)$ with the permutation $(1)=1$, $(2)=3$, $(3)=2$.
Then $x_{(0)}=0$, $x_{(1)}=1$, $x_{(2)}=1$, $x_{(3)}=3$, therefore  $\Psi_{\x}=\{0,2,3\}$ and
% {\bs{nemali by sme v indexe uz pisat namiesto $\x$ $(1,3,1)$, t.j. $\Psi_{(1,3,1)}=\{0,2,3\}$}} and
$$\Ai{1}=\{(1), (2), (3)\}=\{1,2,3\},\quad \Ai{3}=\{(3)\}=\{2\},\quad \Ai{4}=\emptyset. $$
We can see, that the assertion $\mathrm{(C1)}$ of Corollary~\ref{LM_SLM_coincide} is satisfied with 
\begin{center} $G_0=\emptyset$, $G_2=\{3\}$, $G_3=\{2\}$.
\end{center}
Indeed, $\aAi[G_0]{\x}{\mathrm{sum}}=0=x_{(0)}$ and $\mu(G_0^c)=\mu(E_{(1)})$. Further, $\aAi[G_2]{\x}{\mathrm{sum}}=1=x_{(2)}$ and $\mu(G_2^c)=\mu(\{1,2\})=\mu(E_{(3)})$. Finally, $\aAi[G_3]{\x}{\mathrm{sum}}=3=x_{(3)}$ and $\mu(G_3^c)=\mu(\{1,3\})=\mu(E_{(4)})$.
The assertion $\mathrm{(C2)}$ is also satisfied, see the visualisation of generalized survival function via parallel lines in Figure~\ref{LM_eq_SLM_picture}. Discussed survival functions coincide and take the form 
$$\lm=\gsf{\mathrm{sum}}{\x}{\alpha}=\mathbf{1}_{[0,1)}(\alpha)+0{.}5\cdot\mathbf{1}_{[1,3)}(\alpha)$$ for $\alpha\in[0,\infty)$. The plot of (generalized) survival function is in Figure~\ref{LM_eq_SLM_picture}.
\end{example}

\begin{example}\label{example}
Let us consider $\cA^{\mathrm{sum}}=\{\aAi[E]{\cdot}{\mathrm{sum}}: E\in2^{[3]}\}$, and normalized monotone measure $\mu$ on $2^{[3]}$  with the following values:
% Let us consider $n=3$, $\collectionGSF=2^{[3]}$, the input vector $\x=(2,3,4)$ with the permutation  being the identity, $\sA^{\mathrm{sum}}$ and the normalized monotone measure $\mu$ on $2^{[3]}$ with the following values:
\begin{table}[H]
\renewcommand*{\arraystretch}{1.2}
\begin{center}
\begin{tabular}{|c|c|c|c|c|c|c|c|c|}
\hline
$E$ & $\emptyset$ & $\{1\}$ & $\{2\}$ & $\{3\}$ & $\{1,2\}$ & $\{1,3\}$ & $\{2,3\}$ & $\{1,2,3\}$ \\ \hline
$\mu(E)$ & $0$ & $0$ & $0$ & $0.7$ & $0$ & $0.8$ & $0.7$ & $1$\\ \hline
$\aAi[E]{\x}{\mathrm{sum}}$ & $0$  & $2$ & $3$ & $4$ & $5$ & $6$ & $7$ & $9$\\\hline
\end{tabular}
%\caption{The values of conditional aggregation operator $\sA^{\mathrm{sum}}$ from Example~\ref{SC_example}}
%\label{table}
\end{center}
\end{table}
\noindent Further, let us take the input vector $\x=(2,3,4)$ with the permutation  being the identity.
Then survival functions coincide
$$\lm=\gsf{\mathrm{sum}}{\x}{\alpha}=\mathbf{1}_{[0,2)}(\alpha)+0{.}7\cdot\mathbf{1}_{[2,4)}(\alpha).$$
Here, $G_0=\emptyset$, $G_1=\{1\}$, $G_2=\{2\}$, $G_3=\{3\}$ are the only sets that satisfy the equality $\aAi[G_k]{\x}{\mathrm{sum}}=x_{(k)}$ for $k\in\Psi_{\x}=\{0,1,2,3\}$. However, $$0{.}8=\mu(G_2^c)\neq\mu(\Ai{3})=0{.}7.$$
% , which is contradiction with the assumption of Corollary~\ref{LM_SLM_coincide}. 
Thus, a sufficient condition in Corollary~\ref{LM_SLM_coincide} is not a necessary condition. \end{example}

Let us return to Proposition~\ref{LM<SLM}. While $\mathrm{(C2)}$ is the necessary and sufficient condition under which the generalized survival function is greater or equal to the survival function, $\mathrm{(C1)}$ is only sufficient for the reverse inequality. Since this condition seems too strict, % However, the condition $\mathrm{(C1)}$ is not a necessary condition. 
% This is guaranteed by condition (C3). 
% For this reason 
let us define conditions $\mathrm{(C3)}$ and $\mathrm{(C4)}$ as follows:

\begin{enumerate}
   \item[(C3)] For any $k\in\Psi_{\x}$ there exists $F_k\in \collectionGSF$ such that $\aA[F_k]{\x}\leq x_{(k)}$ and $\mu(F_k^c)\leq \mu(\Ai{k+1})$.
   \item[(C4)] For any $k\in\Psi_{\x}$ there exists $F_k\in \collectionGSF$ such that $\aA[F_k]{\x}\leq x_{(k)}$ and $\mu(F_k^c)= \mu(\Ai{k+1})$. 
  \end{enumerate}
  
\begin{figure}
    \begin{center}
    \includegraphics[scale=1.1]{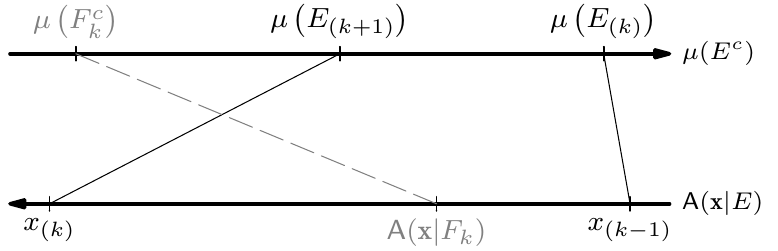}
    \caption{The visualization of condition $\mathrm{(C3)}$ from Proposition~\ref{LM>SLM}}\label{obr_dokaz3}
    \end{center}
\end{figure}

\noindent The visualization of condition $\mathrm{(C3)}$ is drawn in Figure~\ref{obr_dokaz3}. In the following we show that exactly $\mathrm{(C3)}$ improves Proposition~\ref{LM<SLM} $\text{ii)}$. % , $\mathrm{(C3)}$ is sufficient and necessary condition for $\gsf{}{\x}{\alpha}\leq\lm$). 
As a consequence we also get improvement of  Corollary~\ref{LM_SLM_coincide}. Replacing $\mathrm{(C1)}$ with $\mathrm{(C3)}$, we obtain sufficient and necessary condition for equality between survival functions. However, it will turn out that under the $\mathrm{(C2)}$ assumption $\mathrm{(C3)}$ will be reduced to $\mathrm{(C4)}$. 
% \jb{nasledujuca veta nie je sformulovana dobre}Asi nerozumiem, kludne preformuluj.Thus explaining the introduction of this condition on this place.

\begin{proposition}\label{LM>SLM}
Let $\x\in[0,\infty)^{[n]}$, $\mu\in\bM$, and let $\cA$ be FCA.
% The condition 
$\mathrm{(C3)}$ holds if and only if $\gsf{}{\x}{\alpha}\leq\lm$
for any $\alpha\in[0,\infty)$.
\end{proposition}
\proof 
Let us prove the implication $\Rightarrow$. According to Proposition~\ref{vlastnostipi} (iv) let us divide interval $[0,\infty)$ into disjoint sets $$[0,\infty)=\bigcup_{k\in\Psi_{\x}}[x_{(k)},x_{(k+1)}).$$
Let us consider an arbitrary (fixed) $k\in\Psi_{\x}$.
% Let $\alpha\in[x_{(k)},x_{(k+1)})$ for some $k\in\Psi_{\x}$. 
Then by assumptions, there is $F_{k}\in \collectionGSF$ such that $\aA[F_{k}]{\x}\leq x_{(k)}$ and $\mu(F_{k}^c)\leq \mu(\Ai{k+1})$. Thus $\mu(F_k^c)\in\{\mu(E^c):\, \aA[E]{\x}\leq \alpha,\, E\in\collectionGSF \}$ for any 
% {\lh{$\alpha\geq x_{(k)}$. }}
$\alpha\in[x_{(k)},x_{(k+1)})$. 
Hence,
\begin{align*}
\gsf{}{\x}{\alpha}=\min\left\{\mu(E^c):\aA[E]{\x}\leq\alpha, E\in\collectionGSF\right\}\leq\mu(F_k^c)\leq\mu(\Ai{k+1})=\lm\text{}
\end{align*}
for any $\alpha\in[x_{(k)},x_{(k+1)})$.

Let us prove the reverse implication $\Leftarrow$.
Let  $\gsf{}{\x}{\alpha}\leq\lm$ for any $\alpha\in[0,\infty)$. Then from this fact and from~\eqref{PiLM_form} it follows:
$$\gsf{}{\x}{x_{(k)}}\leq\mu(\{\x> x_{(k)}\})=\mu(E_{(k+1)})$$ for any $k\in\Psi_{\x}$. As $\gsf{}{\x}{x_{(k)}}=\min\left\{\mu(E^c): \aA[E]{\x}\leq x_{(k)}, E\in\collectionGSF\right\}$, there exists $F_k\in\collectionGSF$ such that $\aA[F_k]{\x}\leq x_{(k)}$ and $\mu(F_k^c)\leq\mu(\Ai{k+1})$.
%$$\gsf{}{\x}{x_{(k)}}=\min\left\{\mu(E^c): \aA[E]{\x}\leq x_{(k)}, E\in\collectionGSF\right\}\leq\mu(\Ai{k+1})=\mu(\{\x>\alpha\})$$
%for any $k\in\Psi_{\x}$. This implies the fact that  there is $F_k\in\collectionGSF$ such that $\aA[F_k]{\x}\leq x_{(k)}$ and $\mu(F_k^c)\leq\mu(\Ai{k+1})$.
\qed
\bigskip

% Under slight modification of the condition $\mathrm{(C3)}$ we get the result which will be useful in the last proof of this subsection. Let us formulate the further condition as follows:
% \begin{enumerate}[(C4)]
% \item For any $k\in\Psi_{\x}$ there exists $F_k\in \collectionGSF$ such that $\aA[F_k]{\x}\leq x_{(k)}$ and $\mu(F_k^c)= \mu(\Ai{k+1})$. 
% \end{enumerate}

\begin{corollary}\label{corollary_1}
Let $\x\in[0,\infty)^{[n]}$, $\mu\in\bM$, and let $\cA$ be FCA.
\begin{enumerate}[i)]
    \item  If $\mathrm{(C2)}$ holds, then $\mathrm{(C3)}$ is equivalent to $\mathrm{(C4)}$.
    \item $\mathrm{(C2)}$ and $\mathrm{(C3)}$ hold if and only if $\gsf{}{\x}{\alpha}= \lm$
for any $\alpha\in[0,\infty)$.
     \item $\mathrm{(C2)}$ and $\mathrm{(C4)}$ hold if and only if $\gsf{}{\x}{\alpha}= \lm$
for any $\alpha\in[0,\infty)$.
\end{enumerate}
\end{corollary}
\proof It is enough to prove part i), more precisely, the implication $\mathrm{(C3)}\Rightarrow\mathrm{(C4)}$. Let $\mathrm{(C3)}$ is satisfied, we show that $\mu(F_k^c)= \mu(\Ai{k+1})$ holds for any $k\in\Psi_{\x}$. Since for any $F_k\in\collectionGSF$, $k\in\Psi_{\x}$ we have $\aA[x]{F_k}\leq x_{(k)}<x_{(k+1)}$, then from $\mathrm{(C2)}$ we have $\mu(F_k^c)\geq\mu(E_{(k+1)})$. On the other hand, from $\mathrm{(C3)}$ we have $\mu(F_k^c)\leq\mu(E_{(k+1)})$. \qed
\medskip

\begin{remark}\label{remark_i}
At the end of this main part let us remark that some above mentioned results are true also without constructing $\Psi_{\x}$ system. Let us denote:
\begin{enumerate}[$(\widetilde{\mathrm{C}}1)$]
    \item For any $i\in[n]\cup\{0\}$ there exists $G_i\in \collectionGSF$ such that $\aA[G_i]{\x}=x_{(i)}$ and $\mu(F_i^c)=~\mu(\Ai{i+1})$.
    \item {{For any $i\in[n]\cup\{0\}$ and for any  $E\in \collectionGSF$ it holds:}}
    $\aA[E]{\x}<x_{(i+1)}\Rightarrow\mu(E^c)\geq \mu(\Ai{i+1}).$
    \item For any $i\in[n]\cup\{0\}$ there exists $F_i\in \collectionGSF$ such that $\aA[F_i]{\x}\leq x_{(i)}$ and $\mu(F_i^c)\leq \mu(\Ai{i+1})$.
  \end{enumerate}
Then Proposition~\ref{LM<SLM} and Corollary~\ref{corollary_1}~(ii) remain to be true, although, requirements in $(\widetilde{\mathrm{C}}1)$, $(\widetilde{\mathrm{C}}2)$, $(\widetilde{\mathrm{C}}3)$ will be for some  $i\in[n]\cup\{0\}$ redundant \footnote{They will be redundant for $i\in[n]\cup\{0\}$ such that $x_{(i)}=x_{(i+1)}$, compare with the motivation of $\Psi_{\x}$ system introduction.}. On the other hand, Corollary~\ref{corollary_1}~(i),~(iii) need not be satisfied in general.
\bigskip

\noindent{Inequalities:} Let $\x\in[0,\infty)^{[n]}$, $\mu\in\bM$, and let $\cA$ be FCA.
\begin{enumerate}[i)]
\item If $(\widetilde{\mathrm{C}}1)$ holds, then  $\gsf{}{\x}{\alpha}\leq\lm$
for any $\alpha\in[0,\infty)$. %\label{LM>SLM_0}
\item $(\widetilde{\mathrm{C}}2)$ holds if and only if $\lm\leq \gsf{}{\x}{\alpha}$
for any $\alpha\in[0,\infty)$.
\end{enumerate}
\smallskip

\noindent{Sufficient and necessary condition:} Let  $\x\in[0,\infty)^{[n]}$, $\mu\in\bM$, and let $\cA$ be FCA. $(\widetilde{\mathrm{C}}2)$ and $(\widetilde{\mathrm{C}}3)$ hold if and only if $\gsf{}{\x}{\alpha}= \lm$
for any $\alpha\in[0,\infty)$.
\end{remark}

\subsection{Equality of generalized survival function and standard survival function, further results}
\vspace*{12pt}

In this subsection we provide further results on indistinguishability of survival functions. 
Considering the formula of standard survival function~\eqref{PiLM_form} one can observe that the same value of monotone measure may be achieved on several intervals. These intervals can be joined together.
Thus we obtain again a shorter formula of survival function, see Proposition~\ref{LM=SLM_max2} i), which allows us to formulate further results.
% \noindent\textbf{Observation:}
% Let us consider Example~\ref{example} again. Let us
% visualize the standard survival function represented via $\gsf{\textrm{max},\collectionGSF}{\x}{\alpha}$ with $\collectionGSF=\{\Ai{k+1}: k\in\Psi_{\x}=\{0,1,2,3\}\}$ according to~(\ref{PiLM_form}), see Figure~\ref{viz_a_gsf_pred_def_pi_s_hviezdickou}.
% \begin{figure}
% \begin{center}
% \includegraphics[scale=1.1]{obrazky/obrazok_vizualizaciaSF_pred_def_pi_s_hviezdickou.pdf}
% \caption{Generalized survival function and visualization from Example~\ref{example}}  \label{viz_a_gsf_pred_def_pi_s_hviezdickou}
% \end{center}
% \end{figure}
% From the visualization we can see that the line connecting $x_{(2)}$ with $\mu(\Ai{3})$ is not important in calculation of standard survival function, i.e. the same result we get considering collection $\collectionGSF\setminus\Ai{3}$ (equivalently $\Psi_{\x}\setminus\{2\}$).
% Thus, we can reduce the \mbox{system $\Psi_{\x}$.}
% \mbox{}
Let us define
system $\Psi_{\x}^{*}\subseteq\Psi_{\x}$ as follows:
\begin{align}\label{pi*}\Psi_{\x}^*:=\{k\in\Psi_{\x}\setminus\{\min\Psi_{\x}\}:\, \mu(\Ai{j+1})>\mu(\Ai{k+1}), j<k, j\in\Psi_{\x}\}\cup \{\min\Psi_{\x}\}\end{align}
(compare with the definition of system $\Psi_{\x}$ which is analogous, however the main condition is concentrated on components of $\x$ instead of values of $\mu$).
Let us give an example of the $\Psi_{\x}^*$ system calculation considering inputs from Example~\ref{example}. For given input 
$\Psi_{\x}=\{0,1,2,3\}$.
Then by definition of $\Psi_{\x}^*$ we have $\min\Psi_{\x}=0\in\Psi_{\x}^*$.
For $k=1,3$ the inequality $\mu(\Ai{j+1})>\mu(\Ai{k+1})$, $j<k$ holds, however, for $k=2$ we have $\mu(\Ai{2})=\mu(\Ai{3})$.
Thus $2\notin\Psi_{\x}^*$.
In summary, $\Psi_{\x}^*=\{0,1,3\}$. 

For purpose of this subsection for any $k\in\Psi_{\x}^*$  let us denote
\begin{align}\label{l_k}
    l_k:=\max\{j\in\Psi_{\x}:\, \mu(E_{(j+1)})=\mu(E_{(k+1)})\}.
\end{align}

\begin{proposition}\label{vlastnoti_pi_s_hviezdickou}
Let $\x\in[0,\infty)^{[n]}$, $\mu\in\mathbf{M}$.
\begin{enumerate}[i)]
    \item $x_{(\min\Psi_{\x}^*)}=0$.
    \item $\left\{[x_{(k)},x_{(l_k+1)}): k\in\Psi_{\x}^*\right\}$ with $l_k$ given by~\eqref{l_k} and with the convention $x_{(n+1)}=\infty$ is a decomposition of interval $[0,\infty)$ into nonempty pairwise disjoint sets. 
    \item $(\forall k\in\Psi_{\x}^*\setminus\{\min\Psi_{\x}^*\})$ $(\exists r\in\Psi_{\x}^*, r<k)$ $x_{(k)}=x_{(l_{r}+1)}$. 
    Moreover, $\mu(\Ai{k+1})<\mu(\Ai{l_r+1})$.
    % \item  If $j_k=\max\{j\in\Psi_{\x}^*:\,j<k\}$, $k\in\Psi_{\x}^*$, then $l_{j_k}+1=k$ and $\mu(\Ai{k+1})<\mu(\Ai{l_{j_k+1}})$.
   % \item $\bigcup_{k\in\Psi_{\x}^*}[x_{(k)},x_{(l_k)})$ with $l_k=\min\{j\in\Psi_{\x}^*:\,k<j\}$ for any $k\in\Psi_{\x}^*\setminus\{\max\Psi_{\x}^*\}$ with the convention $x_{(l_{\max\Psi_{\x}^*})}$ is a decomposition of interval $[0,\infty)$ into nonempty pairwise disjoint sets
%   \item {\bs{???Platí: $l_k\geq k > l_{r_k}\geq r_k$, $k\in\Psi_{\x}^*$ ???}}\jb{nedavala by som to}
    \item If $\mu$ is such that it is strictly monotone on $\{\Ai{k+1}: k\in\Psi_{\x}\}$, then $\Psi_{\x}=\Psi_{\x}^*$.
\end{enumerate}
\end{proposition}
\proof Part i) follows from Proposition~\ref{vlastnostipi} iii). Since $\Psi_{\x}^*\subseteq\Psi_{\x}$ and $\min\Psi_{\x}\in\Psi_{\x}^*$, then $\min\Psi_{\x}=\min \Psi_{\x}^*$. The proof of ii) follows from Proposition~\ref{vlastnostipi} part iv) and from the fact that each partial interval $[x_{(k)}, x_{(l_k+1)})$,  $k\in\Psi_{\x}^*$ can be rewritten as follows $$[x_{(k)}, x_{(l_k+1)})=\bigcup_{j=k,j\in\Psi_{\x}}^{l_k}[x_{(j)}, x_{(j+1)}).$$ 
The equality $x_{(k)}=x_{(l_{r}+1)}$ in part iii) follows from ii) with $r=\max\{j\in\Psi_{\x}^*: x_{(j)}<x_{(k)}\}$.
Moreover, it holds $\mu(\Ai{l_r+1})=\mu(\Ai{r+1})>\mu(\Ai{k+1})$ where the first equality holds because of~\eqref{l_k}, the second inequality is true due to $r< k$, $r,k\in\Psi_{\x}^*$.
% The inequality $l_{r}<k$ follows from the fact $\mu(\Ai{l_r+1})=\mu(\Ai{r+1})>\mu(\Ai{k+1})$ where the first equality holds because of~\eqref{l_k}, the second inequality is true due to $r< k$, $r,k\in\Psi_{\x}^*$.
Part~iv) follows from~\eqref{pi*}.
\qed

% {\bs{Further, part iii) follows from fact $j_k<k$ and $l_{j_k}=\max\{j\in\Psi_{\x}:\mu(\Ai{j+1})<\mu(\Ai{k+1})\}=k-1$.
% Thus $l_{j_k}+1=k$.
% Then $\mu(\Ai{l_{j_k}+1})=\mu(\Ai{k})>\mu(\Ai{k+1})$ by fact that $k\in\Psi_{\x}^*$.
% }}

\begin{proposition}\label{LM=SLM_max2}
Let $\x\in[0,\infty)^{[n]}$,  $\mu\in\mathbf{M}$. 
\begin{enumerate}[i)]
    \item Then
    \begin{equation}\label{Pi*LM_form}
    \lm=\sum_{k\in\Psi_{\x}^*}\mu\left(\Ai{k+1}\right)\cdot\mathbf{1}_{[x_{(k)},x_{(l_k+1)})}(\alpha)
    \end{equation}
    for any $\alpha\in[0,\infty)$ with $l_k$ given by~\eqref{l_k} and with the convention $x_{(n+1)}=\infty$.
    \item If $\collectionGSF\supseteq\{\Ai{k+1}^c: k\in\Psi_{\x}^*\}$, then 
    $\gsf{\mathrm{max}}{\x}{\alpha}= \lm$
    for any $\alpha\in[0,\infty)$.
   \end{enumerate}
\end{proposition}
\proof Part ii) can be proved analogously as Proposition~\ref{LM=SLM_max} part ii). Part i) follows from the fact that each partial interval $[x_{(k)}, x_{(l_k+1)})$,  $k\in\Psi_{\x}^*$ can be rewritten as follows $[x_{(k)}, x_{(l_k+1)})=\bigcup_{j=k,j\in\Psi_{\x}}^{l_k}[x_{(j)}, x_{(j+1)}).$
    From the formula~\eqref{PiLM_form} and from definition of $l_k$ we get
    $$\lm=\mu(\Ai{j+1})=\mu(\Ai{k+1}).$$
    for any $\alpha\in[x_{(j)}, x_{(j+1)})$.
    \qed
   \bigskip

All results from the previous subsection will also be true under a slight modification of conditions (C1), (C2), (C3) and (C4) as follows:
\begin{enumerate}[(C1$^*$)]
    \item For any $k\in\Psi_{\x}^*$ there exists $G_k\in\collectionGSF$ such that $
    \aA[G_k]{\x}=x_{(k)}$ and  $\mu(G_k^c)=\mu(\Ai{k+1}).$ 
      \item For any $k\in\Psi_{\x}^*$ and for any  $E\in \collectionGSF$ it holds:   $\aA[E]{\x}<x_{(l_k+1)}\Rightarrow\mu(E^c)\geq \mu(\Ai{l_k+1}).$
     \item For any $k\in\Psi_{\x}^*$ there exists $F_k\in \collectionGSF$ such that $\aA[F_k]{\x}\leq x_{(k)}$ and $\mu(F_k^c)\leq \mu(\Ai{k+1})$.
    \item For any $k\in\Psi_{\x}^{*}$ there exists $F_k\in \collectionGSF$ such that $\aA[F_k]{\x}\leq x_{(k)}$ and $\mu(F_k^c)= \mu(\Ai{k+1})$. 
    \end{enumerate}

In the following we summarize all modifications of results from the main part of this section. Since proofs of parts i) -- vii) are based on the same ideas, we omit them. The comparison of these results with those obtained in the main part can be found in Remark~\ref{remark*}.

\begin{proposition}\label{summar}
Let $\x\in[0,\infty)^{[n]}$, $\mu\in\bM$, and let $\cA$ be FCA.
\begin{enumerate}[i)]
\item If $\mathrm{(C1^*)}$ holds, then  $\gsf{}{\x}{\alpha}\leq\lm$
for any $\alpha\in[0,\infty)$.
\item $\mathrm{(C2^*)}$ holds if and only if $\lm\leq \gsf{}{\x}{\alpha}$
for any $\alpha\in[0,\infty)$.
\item If {$\mathrm{(C1^*)}$} and {$\mathrm{(C2^*)}$} are satisfied, then $\gsf{}{\x}{\alpha}= \lm$
for any $\alpha\in[0,\infty)$.
\item $\mathrm{(C3^*)}$ holds  if and only if   $\gsf{}{\x}{\alpha}\leq\lm$
for any $\alpha\in[0,\infty)$.
\item {$\mathrm{(C2^*)}$} and {$\mathrm{(C3^*)}$} hold if and only if  $\gsf{}{\x}{\alpha}= \lm$
for any $\alpha\in[0,\infty)$.
\item If {$\mathrm{(C2^*)}$} holds, then {$\mathrm{(C3^*)}$} is equivalent to {$\mathrm{(C4^*)}$}.
\item {$\mathrm{(C2^*)}$} and {$\mathrm{(C4^*)}$} hold if and only if  $\gsf{}{\x}{\alpha}= \lm$
for any $\alpha\in[0,\infty)$.
\item $\mathrm{(C2)}$ holds if and only if $\mathrm{(C2^*)}$ holds.
\end{enumerate}
\end{proposition}
\proof The implication (C2) $\Rightarrow$ (C2$^*$) of part viii) is clear. We prove the reverse implication. 
Let us consider any set $E\in\collectionGSF$
such that $\aA[E]{\x}<x_{(k+1)}$ for some $k\in\Psi_{\x}$. Let us define 
$$j_k=\min\{j\in\Psi_{\x}: \mu(\Ai{j+1})=\mu(\Ai{k+1})\}.$$
It is easy to see that $j_k\in\Psi_{\x}^*$, $l_{j_k}\geq k\geq j_k$. Moreover,  $\mu(\Ai{l_{j_k}+1})=\mu(\Ai{k+1})=\mu(E_{(j_k+1)})$ and $x_{(k+1)}\leq x_{(l_{j_k}+1)}$. 
Then from (C2$^*$) we have $\mu(E^c)\geq\mu(E_{(l_{j_k}+1)})=\mu(E_{(k+1)})$.
\qed

\begin{remark}\label{remark*}
In comparison with results in the main part of this section, the advantage of previous statements lies in their efficiency for survival functions equality or inequality testing. In particular, Proposition~\ref{summar} vii) requires to hold the same properties as Corollary~\ref{corollary_1} iii), however for a smaller number of sets, $k\in\Psi_{\x}^{*}\subseteq\Psi_{\x}$. On the other hand, the equality (inequality) of survival functions 
% provides 
implies
more information than those included in the Proposition~\ref{summar}, the results are true for any $k\in\Psi_{\x}$ not only for $k\in\Psi_{\x}^*$.  Moreover, system $\Psi_{\x}$ is also easier in definition.
\end{remark}

% Under conditions (C1$^{*}$), (C2$^*$) we are able to formulate sufficient and necessary condition under which (generalized and standard) survival functions coincide without the assumption of strict monotonicity of $\mu$, compare with Corollary~\ref{LM_SLM_coincide}. The assumption of strictly monotonicity of $\mu$ is not needed because of good properties of system $\Psi_{\x}^*$. Since the proof uses very similar ideas like these presented in previous results, we leave it for the Appendix.

We have seen in the main part of this section that $\mathrm{(C1)}$, $\mathrm{(C2)}$ are not necessary for equality between survival functions in general, see Corollary~\ref{LM_SLM_coincide}, Example~\ref{example}.  This result we have improved by replacing $\mathrm{(C1)}$ with $\mathrm{(C4)}$. Also, Corollary~\ref{LM_SLM_coincide} can be improved as it follows. 
%However, from the necessity condition point of view it still makes a sense to formulate necessary conditions of indistinguishability of survival functions with stricter condition than $\mathrm{(C4)}$, e.g. with the condition similar to $\mathrm{(C1)}$, see Remark~\ref{remark_comp}.

\begin{theorem}\label{nutna_postacujuca_2}
Let $\x\in[0,\infty)^{[n]}$, $\mu\in\bM$, and let $\cA$ be FCA. Then the following assertions are equivalent:
\begin{enumerate}[i)]
    \item $\mathrm{(C1^*)}$, $\mathrm{(C2^*)}$ are satisfied.
    \item $\gsf{}{\x}{\alpha}=\lm$ for any $\alpha\in[0,\infty)$.
\end{enumerate}
\end{theorem}
\proof 
The implication $\text{i)}\Rightarrow \text{ii)}$ follows from Proposition~\ref{summar} iii). In order to prove the reverse implication, assume that survival functions are equal. Then (C2$^*$) follows from Proposition~\ref{summar} vii). It is enough to prove (C1$^*$). From Proposition~\ref{summar} vii) we have:
\begin{center}
    For any $k\in\Psi_{\x}^*$ there exists $F_k\in \collectionGSF$ such that $\aA[F_k]{\x}\leq x_{(k)}$ and $\mu(F_k^c)= \mu(\Ai{k+1})$.
\end{center}
We show that $\aA[F_k]{\x}= x_{(k)}$. Indeed, for $k=\min\Psi_{\x}^*$ is the result immediate since 
$0\leq\aA[F_{\min\Psi_{\x}^*}]{\x}\leq~ x_{(\min\Psi_{\x}^*)}=0$, where the last inequality follows from Proposition~\ref{vlastnoti_pi_s_hviezdickou} i). Let $k>\min\Psi_{\x}^*$, $k\in\Psi_{\x}^*$.
From Proposition~\ref{vlastnoti_pi_s_hviezdickou} iii) there exists $r\in\Psi_{\x}^*, r<k$ such that $x_{(l_{r}+1)}=x_{(k)}$ and $\mu(F_k^c)= \mu(\Ai{k+1})< \mu(\Ai{l_{r}+1})$. From contraposition to (C2$^*$) we have $\aA[F_k]{\x}\geq x_{(l_{r}+1)}=x_{(k)}$. 
\qed
\medskip

\begin{corollary}\label{nutna_postacujuca}
Let $\x\in[0,\infty)^{[n]}$, $\mu\in\bM$ {such that it is strictly monotone} on $\{\Ai{k+1}: k\in\Psi_{\x}\}$, and let $\cA$ be FCA. Then the following assertions are equivalent:
\begin{enumerate}[i)]
    \item $\mathrm{(C1)}$, $\mathrm{(C2)}$ are satisfied.
    \item $\gsf{}{\x}{\alpha}=\lm$ for any $\alpha\in[0,\infty)$.
\end{enumerate}
\end{corollary}
\proof 
It follows from Proposition~\ref{vlastnoti_pi_s_hviezdickou} iv) and Theorem~\ref{nutna_postacujuca_2}.
\qed
\medskip

A summary of relationships among some conditions as well as the summary of 
sufficient and necessary conditions under which survival functions coincide or under which they are pointwise comparable with respect to $\leq$, $\geq$ can be found in the Appendix, see Table~\ref{tab:summar}. 

%{\lh{
%\begin{remark}\label{remark_comp}
%Let us compare Corollary~\ref{nutna_postacujuca} with Corollary~\ref{corollary_1} ii), iii). The advantage of stricter necessary condition  $\mathrm{(C1)}$ (in Corollary~\ref{nutna_postacujuca}) in comparison with $\mathrm{(C3)}$, $\mathrm{(C4)}$ (in Corollary~\ref{corollary_1} ii), iii)) is in situations when one wants to show that survival functions do not coincide. In case of Corollary~\ref{nutna_postacujuca} for concluding that functions do not coincide, it is enough to find  $k\in\Psi_{\x}$ such that there does not exist a set $G_k\in\collectionGSF$ with the property $\aA[G_k]{\x}=x_{(k)}$ and $\mu(G_k^c)=\mu(\Ai{k+1})$ (compare with the condition of the Corollary~\ref{corollary_1} iii) $\aA[G_k]{\x}\leq x_{(k)}$ and $\mu(G_k^c)=\mu(\Ai{k+1})$).
%\end{remark}}}

\section{Equality characterization}

% The conditions $\mathrm{(C1)}$, $\mathrm{(C2)}$ stated in Corollary~\ref{LM_SLM_coincide}
Results of the previous section stated conditions depended on FCA $\cA$,
input vector $\x$ and monotone measure $\mu$ to hold equality between survival functions.
Of course, when one changes the monotone measure and the other inputs stay the same, the equality can violate as the following example shows.

\begin{example}\label{ex:main}
Let us consider $\cA^{\mathrm{sum}}=\{\aAi[E]{\cdot}{\mathrm{sum}}: E\in2^{[3]}\}$, and normalized monotone measure $\mu$ on $2^{[3]}$  with the following values:
% Let us consider $n=3$, $\collectionGSF=2^{[3]}$, the input vector $\x=(1,2,1)$,
% the conditional
% aggregation operator $\sA^{\mathrm{sum}}$ with the monotone measures $\mu$ and $\nu$ given as:
\begin{table}[H]
\renewcommand*{\arraystretch}{1.2}
\begin{center}
\begin{tabular}{|c|c|c|c|c|c|c|c|c|}
\hline
$E$ & $\{1,2,3\}$ & $\{2,3\}$ & $\{1,3\}$ & $\{1,2\}$ & $\{3\}$ & $\{2\}$ & $\{1\}$ & $\emptyset$ \\ \hline
$E^c$ & $\emptyset$ & $\{1\}$ & $\{2\}$ & $\{3\}$ & $\{1,2\}$ & $\{1,3\}$ & $\{2,3\}$ & $\{1,2,3\}$ \\ \hline
$\mu(E^c)$ & $0$ & $0$ & $0$ & $0$ & $0$ & $0.5$ & $0.5$ & $1$\\ \hline
$\nu(E^c)$ & $0$ & $0$ & $0$ & $0$ & $0.5$ & $0.5$ & $0.5$ & $1$\\ \hline
$\aAi[E]{\x}{\mathrm{sum}}$ & $4$  & $3$ & $2$ & $3$ & $1$ & $2$ & $1$ & $0$\\\hline
$\aAi[E]{\x}{\mathrm{max}}$ & $2$  & $2$ & $1$ & $2$ & $1$ & $2$ & $1$ & $0$\\\hline
\end{tabular}
\vspace{-18pt}
\end{center}
\end{table}
\noindent Further, let us take the input vector $\x=(1,2,1)$.\noindent Then we can see
$$\gsf{\mathrm{sum}}{\x}{\alpha}=\bin_{[0,1)}(\alpha)=\lm,\quad \alpha\in[0,\infty),$$
but $$\nu_{\cA^{\mathrm{sum}}}(\x,\alpha)=\bin_{[0,1)}(\alpha)+0{.}5\,\bin_{[1,2)}(\alpha)\not=\bin_{[0,1)}(\alpha)=\nu(\{\x>\alpha\}),\quad \alpha\in[0,\infty).$$ %what is keeping with the existence of sets $F$ such that $\aAi[F]{\x}{\mathrm{sum}}\not=\aAi[F]{\x}{\mathrm{max}}$.
\end{example}

In the following we shall find sufficient and necessary conditions on $\cA$ and $\x$ under which survival functions equal for any monotone measure. So, we answer Problem 2, see Theorem~\ref{thm:characterization}, Theorem~\ref{thm:characterization_mon}. In the second step we characterize FCA for which survival functions equal (for any monotone measure and any input vector). We answer Problem 3. 

\begin{theorem}\label{thm:characterization}
Let $\x\in[0,\infty)^{[n]}$, and  $\cA$ be FCA. Then the following assertions are equivalent:
\begin{enumerate}[i)]
    \item $\collectionGSF\supseteq\{\Ai{k+1}^c : k\in\Psi_{\x}\}$ and
    $\aA[E]{\x}=\aAi[E]{\x}{\mathrm{max}}$ for any $E=\Ai{k+1}^c$ with $k\in\Psi_{\x}$, $\aA[E]{\x}\geq\aAi[E]{\x}{\mathrm{max}}$ otherwise.
      \item For each $\mu\in\mathbf{M}$ such that it is strictly monotone on $\{\Ai{k+1} : k\in\Psi_{\x}\}$ it
    %   $2^{[n]}$ 
      holds $\gsf{}{\x}{\alpha}=\lm
    $ for any $\alpha\in[0,\infty)$.
    \item For each $\mu\in\mathbf{M}$ it holds $\gsf{}{\x}{\alpha}=\lm
    $ for any $\alpha\in[0,\infty)$.
\end{enumerate}
\end{theorem}
\proof 
The implication $\text{i)}\Rightarrow \text{iii)}$ we easily prove by Corollary~\ref{LM_SLM_coincide}. Indeed, for any $k\in\Psi_{\x}$ (C1) is satisfied with $G_k=\Ai{k+1}^c$.
% see Remark~\ref{poznamka_max}.
% , since by assumptions we have $$\aA[G_k]{\x}=\aA[\Ai{k+1}^c]{\x}=\aAi[\Ai{k+1}^c]{\x}{\mathrm{max}}=\max\{x_{(1)},\dots,x_{(k)}\}=x_{(k)}.$$
If $\aA[E]{\x}<x_{(k+1)}$, $k\in\Psi_{\x}$ and $E\in\collectionGSF$, then from assumptions we have
$$\aAi[E]{\x}{\mathrm{max}}\leq\aA[E]{\x}<x_{(k+1)}.$$  
%Thus
% following the same argumentation as in Lemma~\ref{LM=SLM_max} 
Then  we get
$E\subseteq \Ai{k+1}^c$, i.e.,
% Thus from the definition of $\sA^{\rm{max}}$ it holds $E\subseteq\{(1),\dots,(k)\}$
% \footnote{Indeed, if $(r)\in E$, $r\geq k+1$, then  $\aAi[E]{\x}{\mathrm{max}}\geq x_{(r)}\geq x_{(k+1)}$ what is a contradiction.},
$E^c\supseteq\Ai{k+1}$ and for each monotone measure $\mu$ we have $\mu(E^c)\geq\mu(\Ai{k+1})$. Thus (C2) is also satisfied.
          
Let us prove the implication $\text{ii)}\Rightarrow \text{i)}$.
Since the assumption holds for any $\mu\colon 2^{[n]}\to[0,\infty)$ such that it is strictly monotone measure on $\{\Ai{k+1}: k\in\Psi_{\x}\}$, it holds for $\mu$ such that it is strictly monotone measure on $2^{[n]}$ (not only on $\{\Ai{k+1}: k\in\Psi_{\x}\}$).
% Let $\gsf{}{\x}{\alpha}=\lm
%   $ for any $\alpha\in[0,\infty)$ and any strictly monotone measure $\mu$.
%   Then 
From Corollary~\ref{nutna_postacujuca} (C1) holds. Moreover, since sets $E_{(k+1)}$ are the only sets with value equal to $\mu(E_{(k+1)})$, we get $G_{k}=\Ai{k+1}^c$ and $\collectionGSF\supseteq\{\Ai{k+1}^c : k\in\Psi_{\x}\}$. So, from (C1) we have 
$$\aA[\Ai{k+1}^c]{\x}=x_{(k)}=\aAi[\Ai{k+1}^c]{\x}{\mathrm{max}}$$ for any $k\in\Psi_{\x}$. Let us prove the second part of~i). Again, if the equality between survival functions holds for any strictly monotone measure $\mu$ on $\{\Ai{k+1}^c : k\in\Psi_{\x}\}$, then it holds for $\mu\colon 2^{[n]}\to[0,\infty)$ being strictly monotone on the above mentioned collection with values: $$\mu(E)=\mu(E_{(k+1)})\,\, \text{for any set}\,\,E\,\,\text{such that}\,\,\aAi[E^c]{\x}{\mathrm{max}}=x_{(k)},\text{ } k\in\Psi_{\x}.\footnote{Given set function is a monotone measure because $\aAi[\emptyset^c]{\x}{\mathrm{max}}=x_{(n)}$, therefore $\mu(E_{(n+1)})=\mu(\emptyset)=0$. Further, let $E_1\subseteq E_2$, i.e.  $E_1^c\supseteq E_2^c$. Then from Proposition~\ref{vlastnostipi} ii) there exist $k_1,k_2\in\Psi_{\x}$ such that $\aAi[E_1^c]{\x}{\mathrm{max}}=x_{(k_1)}$ and $\aAi[E_2^c]{\x}{\mathrm{max}}=x_{(k_2)}$. It is easy to see that $k_1\geq k_2$ and $\mu(E_1)=\mu(\Ai{k_1+1})\leq\mu(\Ai{k_2+1})=\mu(E_2)$.}$$
Let $E\in\collectionGSF$. Then according to Proposition~\ref{vlastnostipi} ii) there exists $k\in\Psi_{\x}\setminus\{0\}$ such that  $$\aAi[E]{\x}{\mathrm{max}}=x_{(k)}.$$
 %=\aAi[E_{k_r+1}^c]{\x}{\mathrm{max}}
Since $\mu$ is strictly monotone on $\{\Ai{k+1}: k\in\Psi_{\x}\}$, then from Proposition~\ref{vlastnoti_pi_s_hviezdickou} iv) we have $\Psi_{\x}=\Psi_{\x}^*$. Further, 
%   and because of definition $l_k$, $k\in\Psi_{\x}$, see TBA we have
%   $\mu(E^c)=\mu(E_{(k_r+1)})<\mu(E_{(l_{k_r}+1)})$.  Because of contraposition of (C2) from Corollary~\ref{nutna_postacujuca} we have $\aA[E]{\x}\geq x_{(l_{j_k}+1)}=x_{(k)}=\aAi[E_{k_r+1}^c]{\x}{\mathrm{max}}=\aAi[E]{\x}{\mathrm{max}}$.
from Proposition~\ref{vlastnoti_pi_s_hviezdickou} i), if $\aAi[E]{\x}{\mathrm{max}}=x_{(\min\Psi_{\x}^*)}=0$ the result is trivial. Let $k>\min\Psi_{\x}^*$. Then from Proposition~\ref{vlastnoti_pi_s_hviezdickou} iii) there exists $r\in\Psi_{\x}^*, r<k$ such that $x_{(l_{r}+1)}=x_{(k)}$ and $\mu(\Ai{k+1})< \mu(\Ai{l_{r}+1})$. Therefore 
$\mu(E^c)=\mu(\Ai{k+1})< \mu(\Ai{l_{r}+1})$ and from contraposition to (C2$^*$) we have $\aA[E]{\x}\geq x_{(l_{r}+1)}=x_{(k)}=\aAi[E]{\x}{\mathrm{max}}$.
 \qed

\begin{remark}From the previous theorem one can see the other sufficient condition under which the standard and generalized survival functions coincide, i.e., the condition i). 
% together with the assumption $\collectionGSF\supseteq\{\Ai{k+1}^c : k\in\Psi_{\x}\}$.
Let us remark that this sufficient condition is more strict than $\mathrm{(C1)}$  and $\mathrm{(C2)}$, i.e., if i) 
% together with the assumption on $\collectionGSF$ 
is satisfied then $\mathrm{(C1)}$, $\mathrm{(C2)}$ are true, however, the reverse implication need not be true in general, see Example~\ref{SC_example}. 
\end{remark}

According to previous result there are vectors for which the equality between survival functions (for any $\mu$) do not lead to $\sA^{\text{max}}$.

%\jb{Ja viem, ze to bude robota navyse, ale nechceme tento priklad vymenit za taky, kde $\cA$ bude napr. FCA taka, ze na mnozinach $E_{(k+1)}^c$ to bude maximum a na ostatnych cosi ine? Len kvoli tomu, ze s takou FCA nikde nepocitame a tu sa to hodi.}
\begin{example}\label{pr_vektor_y}
Let us consider $\cA=\{\aA[E]{\cdot}:E\in2^{[3]}\}$ with conditional aggregation operator from Example~\ref{pr aggr} iii) with $\w=(0.5,0.5,1)$, $\z=(0.5,0.25,1)$. Let us take the input vector $\x=(2,3,4)$.
The values of $\aA[E]{\x}$, $E\in2^{[3]}$ are summarized in following table:
\begin{table}[H]
\renewcommand*{\arraystretch}{1.2}
\begin{center}
\begin{tabular}{|c|c|c|c|c|c|c|c|c|}
\hline
$E$ & $\{1,2,3\}$ & $\{2,3\}$ & $\{1,3\}$ & $\{1,2\}$ & $\{3\}$ & $\{2\}$ & $\{1\}$ & $\emptyset$\\ \hline
$E^c$ & $\emptyset$ & $\{1\}$ & $\{2\}$ & $\{3\}$ & $\{1,2\}$ & $\{1,3\}$ & $\{2,3\}$ & $\{1,2,3\}$ \\ \hline
$\aA[E]{\x}$ & $4$  & $4$ & $4$ & $3$ & $4$ & $6$ & $2$ & $0$\\\hline
\end{tabular}
\end{center}
\end{table}
\vspace{-18pt}
\noindent Then $\Psi_{\x}=\{0,1,2,3\}$ and it holds
\begin{align*}
    \gsf{}{\x}{\alpha}=\lm&=\mu(\{1,2,3\})\cdot\mathbf{1}_{[0,2)}(\alpha)+\mu(\{2,3\})\cdot\mathbf{1}_{[2,3)}(\alpha)+\mu(\{3\})\cdot\mathbf{1}_{[3,4)}(\alpha)\\
    &=\mu(\Ai{1})\cdot\mathbf{1}_{[0,2)}(\alpha)+\mu(\Ai{2})\cdot\mathbf{1}_{[2,3)}(\alpha)+\mu(\Ai{3})\cdot\mathbf{1}_{[3,4)}(\alpha)
\end{align*}
for any $\alpha\in[0,\infty)$ and monotone measure $\mu$.
So, we have shown that there is vector $\x$ and $\cA\neq\cA^\mathrm{max}$ such that $\gsf{}{\x}{\alpha}=\lm$ for any monotone measure $\mu$. Indeed, $6=\aA[\{2\}]{\x}>\aAi[\{2\}]{\x}{\mathrm{max}}=3$ ($\aA[E]{\x}=\aAi[E]{\x}{\mathrm{max}}$ for any $E\in 2^{[3]}\setminus\{\{2\}\}$).
\end{example}

Aggregation functions with the property being bounded from below by the maximum are in the literature called \textit{disjunctive}, see~\cite{GrabischMarichalMesiarPap2009}. 
% For the class of aggregation operators that are {\lh{nondecreasing}} w.r.t sets we get an interesting consequence.
For FCA $\cA$ nondecreasing w.r.t. sets  we get an interesting consequence.

% \begin{lema}
% Let %$\mu$ be a~monotone measure on $[n]$ and  
% $\x\in[0,\infty)^{[n]}$, {\lh{$\cA$ be FCA nondecreasing w.r.t. sets}}. If
%     $\aA[E]{\x}=\aAi[E]{\x}{\mathrm{max}}$ for any $E=\Ai{k+1}^c$ with $k\in\Psi_{\x}$, otherwise $\aA[E]{\x}\geq\aAi[E]{\x}{\mathrm{max}}$, then 
%     % the assertion (i) in Theorem~\ref{thm:characterization} is equivalent to the assertion
%     $\aA[E]{\x}=\aAi[E]{\x}{\mathrm{max}}$ for any set $E\in\collectionGSF$.
% \end{lema}

\begin{lema}\label{lema_char}
Let $\cA=\{\aA[E]{\cdot}:E\in 2^{[n]}\}$ be FCA nondecreasing w.r.t. sets. If for $\x\in[0,\infty)^{[n]}$ it holds that $\aA[E]{\x}=\aAi[E]{\x}{\mathrm{max}}$ for any $E=\Ai{k+1}^c$ with $k\in\Psi_{\x}$ and $\aA[E]{\x}\geq\aAi[E]{\x}{\mathrm{max}}$ otherwise, then $\aA[E]{\x}=\aAi[E]{\x}{\mathrm{max}}$ for any $E\in2^{[n]}$.\end{lema}

\proof
Let us consider an arbitrary set $E\in\collectionGSF$ and let us denote $\aAi[E]{\x}{\mathrm{max}}:=x_{s}$. Then according to Proposition~\ref{vlastnostipi} there exists $k_s\in\Psi_{\x}$ such that $x_s=x_{(k_s)}$. Then $E\subseteq\Ai{k_s+1}^c$. From above mentioned and from Theorem~\ref{thm:characterization} we have
$$x_{(k_s)}=\aAi[\Ai{k_s+1}^c]{\x}{\mathrm{max}}=\aA[\Ai{k_s+1}^c]{\x}\geq\aA[E]{\x}\geq\aAi[E]{\x}{\mathrm{max}}= x_{(k_s)}.$$
% Let us consider any set $E\in\collectionGSF$. Then any $j\in E$ can be expressed via permutation as $j=(r_j)$. Let us denote $$r_{j}^*=\max_{j\in E}r_j.$$
% Then $\aAi[E]{\x}{\mathrm{max}}=x_{(r_{j}^*)}$ and according to Remark~\ref{remark} there exists $r_k\in\Psi_{\x}$ such that $x_{(r_j^*)}=x_{(r_{k})}$, $r_{j}^*\leq r_{k}$. Then it is clear that $E\subseteq\Ai{r_{j}^*+1}^c\subseteq \Ai{r_{k}+1}^c$. Then we get \begin{align}\label{1}\aAi[E_{(r_k+1)}^c]{\x}{\mathrm{max}}=x_{(r_k)}=x_{(r_j^*)}=\aAi[E]{\x}{\mathrm{max}}\leq\aA[E]{\x}\leq\aA[E_{(r_k+1)}^c]{\x},\end{align} where the last inequality is satisfied because of monotonicity w.r.t. sets of conditional aggregation operator. Since $\aA[E_{(r_k+1)}^c]{\x}=\aAi[E_{(r_k+1)}^c]{\x}{\mathrm{max}}$, from~(\ref{1}) we obtain the required equality.
\qed

\begin{theorem}\label{thm:characterization_mon}
Let %$\mu$ be a~monotone measure on $[n]$ and  
$\x\in[0,\infty)^{[n]}$, and $\cA$ be FCA nondecreasing w.r.t. sets. Then the following assertions are equivalent:
\begin{enumerate}[i)]
    \item $\collectionGSF\supseteq\{\Ai{k+1}^c : k\in\Psi_{\x}\}$ and $\aA[E]{\x}=\aAi[E]{\x}{\mathrm{max}}$ for any set $E\in\collectionGSF$.
    \item For each $\mu\in\mathbf{M}$ it holds $\gsf{}{\x}{\alpha}=\lm
    $ for any $\alpha\in[0,\infty)$.
\end{enumerate}
%In such case, the conditions $\rm{(C1)}$, $\rm{(C2)}$ hold with $B_i=A_{(i+1)}$  for any $i=0,1,2,\dots, n$.
\end{theorem}

% \jb{\begin{theorem}Let $\cA=\{\aA[E]{\x}:E\in 2^{[n]}\}$ be FCA nondecreasing w.r.t. sets. Then the following assertions are equivalent:
% \begin{enumerate}[i)]
%     \item For $\x\in[0,\infty)^{[n]}$ it holds that $\aA[E]{\x}=\aAi[E]{\x}{\mathrm{max}}$ for any set $E\in2^{[n]}$.
%     \item For each monotone measure $\mu$ it holds $\gsf{}{\x}{\alpha}=\lm
%     $ for any $\alpha\in[0,\infty)$.\end{enumerate}
% \end{theorem}}
\proof 
The implication $\text{i)}\Rightarrow \text{ii)}$ follows from Proposition~\ref{LM=SLM_max} ii). The reverse implication follows from Theorem~\ref{thm:characterization} and Lemma~\ref{lema_char}.
\qed
\bigskip

Let us return to Example~\ref{pr_vektor_y}. We have shown that for the input vector $\x=(2,3,4)$ with $\cA$ given in example $\gsf{}{\x}{\alpha}=\lm$ for any $\mu$. However, for another vector, let us take $\y=(2,5,4)$ the equality can violate:
\begin{align*}
    \gsf{}{\y}{\alpha}&=\mu(\{1,2,3\})\cdot\mathbf{1}_{[0,2)}(\alpha)+\mu(\{2,3\})\cdot\mathbf{1}_{[2,4)}(\alpha),\\
    \mu(\{\y>\alpha\})&=\mu(\{1,2,3\})\cdot\mathbf{1}_{[0,2)}(\alpha)+\mu(\{2,3\})\cdot\mathbf{1}_{[2,4)}(\alpha)+\mu(\{2\})\cdot\mathbf{1}_{[4,5)}(\alpha),
\end{align*}
i.e. $\gsf{}{\y}{\alpha}=\mu(\{\y>\alpha\})$ does not hold for any $\mu$.
In the following we shall ask for FCA $\cA$ for which the equality holds for any $\mu$ and for any $\x$. As a last thus we solve Problem 3.

 \begin{theorem}\label{thm:characterization2}
Let $\cA$ be FCA. The following assertions are equivalent:
\begin{enumerate}[i)]
    \item $\mathcal{A}=\{\aAi[E]{\cdot}{\max}:\,\, E\in2^{[n]}\}$. 
    % $\aA[E]{\x}=\aAi[E]{\x}{\mathrm{max}}$ for any $\x\in[0,\infty)^{[n]}$ and for any $E\in\collectionGSF=2^{[n]}$. 
    \item  For each $\mu\in\mathbf{M}$, and for each $\x\in[0,\infty)^{[n]}$ it holds $\gsf{}{\x}{\alpha}=\lm$ for any $\alpha\in[0,\infty)$.
\end{enumerate}
\end{theorem}

\proof The implication $\text{i)}\Rightarrow \text{ii)}$ is immediate. We prove $\text{ii)}\Rightarrow \text{i)}$. Since the equality holds for any $\x$, according to Theorem~\ref{thm:characterization} we get $$\collectionGSF=\bigcup_{\x\in[0,\infty)^{[n]}}\collectionGSF^{\Psi_{\x}-\text{chain}}=2^{[n]}
$$ with $\collectionGSF^{\Psi_{\x}-\text{chain}}:=\{\Ai{k+1}^c: k\in\Psi_{\x}\}.$

Let $\x\in[0,\infty)^{[n]}$ be an arbitrary fixed vector. 
% Let us denote $$\collectionGSF_{{\x}}:=\{\Ai{k+1}^c: k\in\Psi_{\x}\}.$$ 
From Theorem~\ref{thm:characterization} we have $\aA[E]{\x}=\aAi[E]{\x}{\mathrm{max}}$ for any $E\in\collectionGSF^{\Psi_{\x}-\text{chain}}$  and $\aA[E]{\x}\geq\aAi[E]{\x}{\mathrm{max}}$ for any $E\in2^{[n]}\setminus\collectionGSF^{\Psi_{\x}-\text{chain}}$. However, we show that $\aA[E]{\x}=\aAi[E]{\x}{\mathrm{max}}$ for any $E\in2^{[n]}$. 
Let us consider an arbitrary fixed $E\in2^{[n]}\setminus\collectionGSF^{\Psi_{\x}-\text{chain}}$ and vector $$\widehat{\x}=\x\bin_E+a\bin_{E^c},\,\, a>\max_{i\in E}x_i.$$
The set $[n]\in\collectionGSF^{\Psi_{\x}-\text{chain}}$ by definition of $\Psi_{\x}$, therefore $\widehat{\x}\neq \x$. %$E^c\neq\emptyset$. 
Moreover, there exists permutation $(\cdot)$ such that $0=\widehat{x}_{(0)}\leq\widehat{x}_{(1)}\leq\widehat{x}_{(2)}\leq\dots=\widehat{x}_{(\widehat{k})}<\widehat{x}_{(\widehat{k}+1)} =\dots= \widehat{x}_{(n)}=a$ with $\widehat{k}=|E|$. Therefore $\widehat{k}\in\Psi_{\widehat{\x}}$, and $E=\{(1),\dots,(\widehat{k})\}=E^c_{(\widehat{k}+1)}\in\collectionGSF^{\Psi_{\widehat{\x}}-\text{chain}}$.
Finally, from Theorem~\ref{thm:characterization}, and because of the property $\aA[E]{\y}=\aA[E]{\y\bin_E}$ for any $\y\in[0,\infty)^{[n]}$, see~\cite{BoczekHalcinovaHutnikKaluszka2020}, we have: 
\begin{align*}\aA[E]{\x}&=\aA[E]{\x\bin_E}=\aA[E]{\widehat{\x}\bin_E}=\aA[E]{\widehat{\x}}=\aAi[E]{\widehat{\x}}{\mathrm{max}}=\aAi[E]{\widehat{\x}\bin_E}{\mathrm{max}}\\&=\aAi[E]{\x\bin_E}{\mathrm{max}}=\aAi[E]{\x}{\mathrm{max}}.\end{align*}
\qed
% The previous Theorem may be seen as the characterization of FCA for which the survival functions equality holds (for any $\mu$, for any $\x$).

\section{Conclusion}

In this paper we have solved three problems dealing with the question of equality between 
% discussed several conditions to obtain the equality of 
the survival function and the generalized survival function based on conditional aggregation operators introduced originally in~\cite{BoczekHalcinovaHutnikKaluszka2020} (the generalization of concepts of papers~\cite{DoThiele2015},~\cite{HalcinovaHutnikKiselakSupina2019}). We have restricted ourselves to discrete settings. The most interesting results are
% sufficient and necessary conditions under which survival functions coincide, see
Corollary~\ref{LM_SLM_coincide}, 
Corollary~\ref{corollary_1}, Proposition~\ref{summar} and Theorem~\ref{nutna_postacujuca_2} (solutions of Problem 1), Theorem~\ref{thm:characterization} and Theorem~\ref{thm:characterization_mon} (solution of Problem 2). Results were derived from
%based on 
the well-known formula of the standard survival function with a permutation $(\cdot)$ playing a crucial role.
% A permutation $(\cdot)$, with the property of nondecreasing arrangement of the components of input vector, plays the main role in our analysis.
% \jb{In contrast to~\cite{BoczekHalcinovaHutnikKaluszka2020}, we have restricted ourselves to the conditional aggregation of finite number of real inputs only, i.e., we were working on a finite set $[n]:=\{1,2,\dots, n\},\,\, n\geq 2$ with discrete topology. Apart from some helpful reductions, this simplified model allowed us to introduce a system $\Psi_{\x}:=\{i\in[n-1]\cup\{0\}: x_{(i)}<x_{(i+1)}\}\cup\{n\}$, where $(\cdot)$ is a fixed permutation of input vector $\x$ and consequently determine the specific sets $E_{(k+1)}^c$, $k\in\Psi_{\x}$, which play an important role in results throughout this paper.}}
As the main result, we have determined
the family of conditional aggregation operators with respect to which the novel survival function is identical to the standard survival function regardless of the monotone measure and input vector, see Theorem~\ref{thm:characterization2}. 
% As a main result, we have exactly determined
% a class of conditional aggregation operators, which the novel survival function based on, is identical to the standard survival function regardless of the monotone measure. 
% Namely, the values of such conditional aggregation operators have to be equal to values of maximum operator $\sA^{\mathrm{max}}$ on sets $E_{(k+1)}^c$, $k\in\Psi_{\x}$. On sets which are not of the form $E_{(k+1)}^c$, $k\in\Psi_{\x}$, the aggregation operator values need not be smaller than $\sA^{\mathrm{max}}$-values. 

We expect the future extension of our results into the area of integrals introduced with respect to novel survival functions, see~\cite[Definition 5.1]{BoczekHalcinovaHutnikKaluszka2020}. The relationship of studied survival functions (in the sense of equalities or inequalities) determines also the relationship of corresponding integrals (based on standard and generalized survival function). The interesting question for the future work is: Is $\cA^{\mathrm{sup}}$ family of conditional operators also the only one that generates the standard survival function in case of arbitrary basic set $X$ instead of $[n]$, i.e., is it true that

\begin{center}
$\gsf{}{f}{\alpha}=\mu(\{f>\alpha\})$, $\alpha\in[0,\infty)$ for
any $\mu$ and for any $f$ if and only if $\cA=\cA^{\mathrm{sup}}$?
\end{center}

Up to now there are not known any other families except of $\cA^{\mathrm{sup}}$ generating generalized survival function indistinguishable from survival function (for any $\mu$, $f$). 
We believe that new results will be beneficial in some applications, e.g. in the theory of decision making.
The equality between survival functions of a given alternative $\x$ means that the overall score of it with respect to the Choquet integral and the $\cA$-Choquet integral is the same. Also, immediately with decision making application the question of $(\mu,\cA)$-indistinguishability arises, i.e. under which condition on $\mu$, $\cA$ it holds $\gsf{}{\x}{\alpha}=\gsf{}{\y}{\alpha}$ for $\x,\y\in[0,\infty)^{[n]}$. Then alternatives $\x,\y$ will be $\cA$-Choquet integral indistinguishable, i.e. they achieve the same overall score.

\section*{Appendix}
In this subsection we summarize all sufficient and necessary conditions for equality or inequality between survival functions, see Table~\ref{tab:summar}.
% In this subsection we provide the proof of Theorem~\ref{nutna_postacujuca_2} in detail. We also summarize all sufficient and necessary conditions for equality or inequality between survival functions, see Table~\ref{tab:summar}.
\bigskip

\begin{table}[H]
\renewcommand*{\arraystretch}{1.2}
    \centering
    \small
     \begin{tabular}{|c|c|c|l|l|}
    % $\collectionGSF\supseteq\{E_{(k+1)}^c:k\in\Psi_{\x}\}$ & $\Rightarrow$ &
    % $\gsf{\mathrm{max}, {\collectionGSF}}{\x}{\alpha}= \lm$
    % & \\
    \hline
    $(\widetilde{\mathrm{C}}1)$ and $(\widetilde{\mathrm{C}}2)$ & $\Rightarrow$ & $\gsf{}{\x}{\alpha}= \lm$ & & Rem.~\ref{remark_i} \\
    \hline
   $(\widetilde{\mathrm{C}}2)$ and $(\widetilde{\mathrm{C}}3)$ & $\Leftrightarrow$ & $\gsf{}{\x}{\alpha}= \lm$ & & Rem.~\ref{remark_i}\\
    \hline
    $\mathrm{(C1)}$ and $\mathrm{(C2)}$ & $\Rightarrow$ & $\gsf{}{\x}{\alpha}= \lm$ & & Cor.~\ref{LM_SLM_coincide} \\
    \hline
    $\mathrm{(C2)}$ and $\mathrm{(C3)}$ & $\Leftrightarrow$ & $\gsf{}{\x}{\alpha}= \lm$ & & Cor.~\ref{corollary_1}\\
    \hline
    $\mathrm{(C2)}$ and $\mathrm{(C4)}$ & $\Leftrightarrow$ & $\gsf{}{\x}{\alpha}= \lm$ & & Cor.~\ref{corollary_1}\\
      \hline
    \multirow{2}{*}{
    \shortstack[c]{
     $\mathrm{(C1)}$ and $\mathrm{(C2)}$
    }
    }
    &   \multirow{2}{*}{
    \shortstack[c]{$\Rightarrow$}} &   \multirow{2}{*}{
    \shortstack[c]{$\gsf{}{\x}{\alpha}= \lm$ }}& & \multirow{2}{*}{
    \shortstack[l]{\hspace{-3pt}Cor.~\ref{nutna_postacujuca}}}\\
           & & & \multirow{-2}{*}{%
    \shortstack[l]{$\mu\colon2^{[n]}\to[0,\infty)$ is strictly\\ monotone on $\{\Ai{k+1}: k\in\Psi_{\x}\}$}} & \\
    \hline
    $\mathrm{(C2^*)}$ and $\mathrm{(C3^*)}$ & $\Leftrightarrow$ & $\gsf{}{\x}{\alpha}= \lm$ & & Prop.~\ref{summar}\\
    \hline
    $\mathrm{(C2^*)}$ and $\mathrm{(C4^*)}$ & $\Leftrightarrow$ & $\gsf{}{\x}{\alpha}= \lm$ & & Prop.~\ref{summar}\\
    \hline
    $\mathrm{(C1^*)}$ and $\mathrm{(C2^*)}$ & $\Leftrightarrow$ & $\gsf{}{\x}{\alpha}= \lm$ & & Th.~\ref{nutna_postacujuca_2}\\
    \hline
    $\mathrm{(C1)}$ & $\Rightarrow$ & $\gsf{}{\x}{\alpha}\leq \lm$ & & Prop.~\ref{LM<SLM}\\
    \hline
    $\mathrm{(C1^*)}$ & $\Rightarrow$ & $\gsf{}{\x}{\alpha}\leq \lm$ & & Prop.~\ref{summar}\\
    \hline
    $\mathrm{(C3)}$ & $\Leftrightarrow$ & $\gsf{}{\x}{\alpha}\leq \lm$ & & Prop.~\ref{LM>SLM}\\
    \hline
    $\mathrm{(C3^*)}$ & $\Leftrightarrow$ & $\gsf{}{\x}{\alpha}\leq \lm$ & & Prop.~\ref{summar}\\
    \hline
    $\mathrm{(C2)}$ & $\Leftrightarrow$ & $\gsf{}{\x}{\alpha}\geq \lm$ & & Prop.~\ref{LM<SLM}\\
    \hline
    $\mathrm{(C2^*)}$ & $\Leftrightarrow$ & $\gsf{}{\x}{\alpha}\geq \lm$ & & Prop.~\ref{summar}\\
    % \hline
    % %
    % \multirow{2}{*}{
    % \shortstack[c]{
    % $\displaystyle\aA[\{i\}]{\x}\geq x_{i}$\\
    % $\forall i\in[n]$
    % }
    % }
    % & $\Rightarrow$ & $\gsf{}{\x}{\alpha}\geq \lm$ & & Prop.~\ref{LM<=SLM}\\
    % & & & & \\
    % & & & \multirow{-3}{*}{%
    % \shortstack[l]{$\cA$ is monotone w.r.t. sets,\\
    % $\collectionGSF\supseteq\{\Ai{k+1}^c: k\in\Psi_{\x}\}$\\$\,\,\,\,\,\,\,\,\,\,\,\,\,\,\cup\{\{i\}: i\in[n]\}$}} & \\
    % \hline
    % %
    % \multirow{2}{*}{
    % \shortstack[c]{
    % $\displaystyle\aA[\{\Ai{k+1}^c\}]{\x}\leq x_{i}$\\
    % $\forall k\in\Psi_{\x}$
    % %  $\forall k\in\Psi_{\x}\setminus\{0\}$
    % }
    % }
    % & $\Rightarrow$ & $\gsf{}{\x}{\alpha}\leq \lm$ & & Prop.~\ref{veta_p}\\
    % & & & \multirow{-2}{*}{%
    % \shortstack[l]{$\cA$ is monotone w.r.t. sets,\\ $\collectionGSF\supseteq\{\Ai{k+1}^c: k\in\Psi_{\x}\}$}} & \\
    \hline
        \end{tabular}
    \caption{Sufficient and necessary conditions for  pointwise comparison of survival functions}
    \label{tab:summar}
\end{table}

From the Table~\ref{tab:summar} the following relationships between conditions (C1), (C2), (C3), (C4) and its $^*$ versions hold.

\begin{corollary}\label{corollary_2}
Let $\x\in[0,\infty)^{[n]}$, $\mu\in\bM$, and let $\cA$ be FCA. Then it holds:
\begin{align*}
\big(\mathrm{(C1)} \wedge \mathrm{(C2)}\big) & \Rightarrow 
\big(\mathrm{(C1^*)} \wedge \mathrm{(C2^*)}\big) \Leftrightarrow
\big(\mathrm{(C2)} \wedge \mathrm{(C3)}\big) \Leftrightarrow
\big(\mathrm{(C2)} \wedge \mathrm{(C4)}\big)
\Leftrightarrow
\big((\widetilde{\mathrm{C}}2) \wedge (\widetilde{\mathrm{C}}3)\big)
\\&\Leftrightarrow
\big(\mathrm{(C2^*)} \wedge \mathrm{(C3^*)}\big)  \Leftrightarrow
\big(\mathrm{(C2^*)} \wedge \mathrm{(C4^*)}\big)  \Leftrightarrow
\big(\mathrm{(C1^*)} \wedge \mathrm{(C2)}\big).
\end{align*}

% Then the following conditions are equivalent
% \begin{enumerate}[i)]
%     \item  $(C2)$ and $(C3)$ hold,
%      \item  $(C2)$ and $(C4)$ hold,
%       \item  $(C2^*)$ and $(C3^*)$ hold,
%      \item  $(C2^*)$ and $(C4^*)$ hold,
%     \item $(C1^*)$ and $(C2^*)$ hold,
%      \item $(C1^*)$ and $(C2)$ hold,
% %     \item $\gsf{}{\x}{\alpha}= \lm$
% % for any $\alpha\in[0,\infty)$.
% \end{enumerate}
\end{corollary}

\begin{corollary}\label{corollary_3}
Let $\x\in[0,\infty)^{[n]}$, $\mu\in\bM$, and let $\cA$ be FCA. If $\mathrm{(C2^*)}$ holds, then
$$\mathrm{(C1^*)}\Leftrightarrow\mathrm{(C3^*)}\Leftrightarrow\mathrm{(C4^*)}.$$
\end{corollary}

\vspace{1cm}

%\section*{Acknowledgement}
%{\color{red} Authors would like to express their sincere thanks and gratitude to Associate Editor and  anonymous reviewers for their thoughtful suggestions toward the improvement of the paper.}
%This work was supported by the Slovak Research and Development Agency under the  contract  No.  APVV-16-0337. The work is also cofinanced by bilateral call Slovak-Poland grant scheme No. SK-PL-18-0032.

\section*{References}
% \bibliographystyle{siam}
% \bibliography{references}

\begin{thebibliography}{10}

\bibitem{BoczekHalcinovaHutnikKaluszka2020}
{\sc M.~Boczek, L.~Hal\v{c}inov\'{a}, O.~Hutn\'{i}k, and M.~Kaluszka}, {\em
  Novel survival functions based on conditional aggregation operators}, Inform.
  Sciences,  (https://doi.org/10.1016/j.ins.2020.12.049).

\bibitem{BoczekHovanaHutnikKaluszka2020b}
{\sc M.~Boczek, A.~Hovana, O.~Hutn\'ik, and M.~Kaluszka}, {\em New monotone
  measure-based integrals inspired by scientific impact problem}, European
  Journal of Operational Research, 290 (2021), pp.~346--357.

\bibitem{BorzovaHalcinovaSupina2018}
{\sc J.~Borzov{\'a}, L.~Hal{\v{c}}inov{\'a}, and J.~{\v{S}}upina}, {\em
  Size-based super level measures on discrete space}, Medina J. et al. (eds)
  Information Processing and Management of Uncertainty in Knowledge-Based
  Systems. Theory and Foundations. IPMU 2018. Communications in Computer and
  Information Science,  (2018), pp.~219--230.

\bibitem{BorzovaHalcinovaSupina2021}
{\sc J.~Borzov{\'a}, L.~Hal{\v{c}}inov{\'a}, and J.~\v{S}upina}, {\em
  Conditional aggregation-based \uppercase{C}hoquet integral on discrete
  space},  (submitted).

\bibitem{BorzovaHalcinovaHutnik2020}
{\sc J.~Borzov\'{a}, L.~Hal\v{c}inov\'{a}, and O.~Hutn\'{i}k}, {\em The
  smallest semicopula-based universal integrals: Remarks and improvements},
  Fuzzy Sets and Systems, 393 (2020), pp.~29--52.
\newblock Copulas and Related Topics.

\bibitem{CalvoKolesarovaKomornikovaMesiar2002}
{\sc T.~Calvo, A.~Koles{\'a}rov{\'a}, M.~Komorn{\'\i}kov{\'a}, and R.~Mesiar},
  {\em Aggregation operators: Properties, classes and construction methods,
  aggregation operators. new trends and applications}, Physica, Heidelberg,
  (2002), pp.~3--104.

\bibitem{ChenMesiarLiStupnanova2017}
{\sc T.~Chen, R.~Mesiar, J.~Li, and A.~Stup{\v{n}}anov{\'a}}, {\em Possibility
  and necessity measures and integral equivalence}, International Journal of
  Approximate Reasoning, 86 (2017), pp.~62--72.

\bibitem{DoThiele2015}
{\sc Y.~Do and C.~Thiele}, {\em $l^p$ theory for outer measures and two themes
  of lennart carleson united}, Bulletin of the American Mathematical Society,
  52 (2015), pp.~249--296.

\bibitem{DuranteSempi2015}
{\sc F.~Durante and C.~Sempi}, {\em Principles of Copula Theory}, CRC Press, 07
  2015.

\bibitem{grabisch2016set}
{\sc M.~Grabisch}, {\em Set Functions, Games and Capacities in Decision
  Making}, Theory and Decision Library C, Springer International Publishing,
  2016.

\bibitem{GrabischMarichalMesiarPap2009}
{\sc M.~Grabisch, J.~Marichal, R.~Mesiar, and E.~Pap}, {\em Aggregation
  Functions}, Encyclopedia of Mathematics and its Applications, Cambridge
  University Press, 2009.

\bibitem{Halcinova2017}
{\sc L.~Hal{\v{c}}inov{\'a}}, {\em Sizes, super level measures and integrals},
  in Aggregation Functions in Theory and in Practice, 9th International Summer
  School on Aggregation Functions, Sk{\"{o}}vde, Sweden, 19-22 June 2017,
  V.~Torra, R.~Mesiar, and B.~D. Baets, eds., vol.~581 of Advances in
  Intelligent Systems and Computing, Springer, 2017, pp.~181--188.

\bibitem{HalcinovaHutnikKiselakSupina2019}
{\sc L.~Hal{\v{c}}inov{\'a}, O.~Hutn{\'i}k, J.~Kise{\v{l}}{\'a}k, and
  J.~{\v{S}}upina}, {\em Beyond the scope of super level measures}, Fuzzy Sets
  and Systems, 364 (2019), pp.~36--63.

\bibitem{KlementMesiarPap2010}
{\sc E.~P. Klement, R.~Mesiar, and E.~Pap}, {\em {A~universal integral as
  common frame for Choquet and Sugeno integral}}, IEEE Transactions on Fuzzy
  Systems, 18 (2010), pp.~178--187.

\bibitem{Weber1986}
{\sc S.~Weber}, {\em Two integrals and some modified versions - critical
  remarks}, Fuzzy Sets and Systems, 20 (1986), pp.~97--105.

\end{thebibliography}

\end{document}